# RANK-BASED OPTIMAL TESTS OF THE ADEQUACY OF AN ELLIPTIC VARMA MODEL

By Marc Hallin and Davy Paindaveine[1]

*Université Libre de Bruxelles*

We are deriving optimal rank-based tests for the adequacy of a vector autoregressive-moving average (VARMA) model with elliptically contoured innovation density. These tests are based on the ranks of *pseudo-Mahalanobis distances* and on normed residuals computed from Tyler's [*Ann. Statist.* **15** (1987) 234–251] scatter matrix; they generalize the univariate signed rank procedures proposed by Hallin and Puri [*J. Multivariate Anal.* **39** (1991) 1–29]. Two types of optimality properties are considered, both in the local and asymptotic sense, a la Le Cam: (a) (fixed-score procedures) local asymptotic minimaxity at selected radial densities, and (b) (estimated-score procedures) local asymptotic minimaxity uniform over a class $\mathcal{F}$ of radial densities. Contrary to their classical counterparts, based on cross-covariance matrices, these tests remain valid under arbitrary elliptically symmetric innovation densities, including those with infinite variance and heavy-tails. We show that the AREs of our fixed-score procedures, with respect to traditional (Gaussian) methods, are the same as for the tests of randomness proposed in Hallin and Paindaveine [*Bernoulli* **8** (2002b) 787–815]. The multivariate serial extensions of the classical Chernoff–Savage and Hodges–Lehmann results obtained there thus also hold here; in particular, the van der Waerden versions of our tests are uniformly more powerful than those based on cross-covariances. As for our estimated-score procedures, they are fully adaptive, hence, uniformly optimal over the class of innovation densities satisfying the required technical assumptions.

## 1. Introduction.

1.1. *Multivariate signs and ranks.* Much attention has been given recently to the development of invariant, distribution-free and robust methods

Received April 2002; revised May 2003.
[1]Supported in part by a P.A.I. contract of the Belgian Federal Government and an Action de Recherche Concertée of the Communauté française de Belgique.
*AMS 2000 subject classifications.* 62G10, 62M10.
*Key words and phrases.* VARMA models, elliptical symmetry, multivariate ranks, local asymptotic normality.

in the context of multivariate analysis. Whereas such concepts as medians, quantiles, ranks or signs have been present in the classical toolkit of univariate statistical inference for about half a century, the emergence of their multivariate counterparts has been considerably slower.

A fairly complete theory of rank and sign methods for multivariate analysis was elaborated in the sixties, culminating in the monograph by Puri and Sen (1971). This theory however suffers the major weakness of being based on componentwise definitions of ranks and signs, yielding procedures that heavily depend on the choice of a particular coordinate system. It took about twenty years to see the emergence of a systematic development of coordinate-free, affine-invariant competitors to these componentwise sign and rank methods.

This development, initiated in the late eighties, essentially expanded along two distinct lines of research. The first one, based on the so-called *Oja signs and ranks*, is due to Oja, Hettmansperger and their collaborators [Möttönen and Oja (1995), Möttönen, Oja and Tienari (1997), Möttönen, Hettmansperger, Oja and Tienari (1998), Hettmansperger, Nyblom and Oja (1994) and Hettmansperger, Möttönen and Oja (1997); see Oja (1999) for a review]. The second one is associated with Randles' concept of interdirections, and was developed by Randles and his coauthors [Randles (1989), Peters and Randles (1990), Jan and Randles (1994) and Um and Randles (1998)]. For both groups of methods, only location problems (one- and two-sample problems, analysis of variance, ...) were considered, and optimality issues were not investigated. This problem of optimality has been addressed for the first time in Hallin and Paindaveine (2002a) who, still for the location problem, are constructing fully affine-invariant methods based on Randles' interdirections and the so-called pseudo-Mahalanobis ranks that are also fully efficient (in the Le Cam sense) for the multivariate one-sample location model. Invariance and robustness on one side, efficiency on the other, thus, should not be perceived as totally irreconcilable objectives.

The case of multivariate time series problems, in this respect, is much worse, despite the recognized need for invariant, distribution-free and robust methods in the area. The univariate context has been systematically explored, with a series of papers by Hallin, Ingenbleek and Puri (1985), Hallin and Mélard (1988) and Hallin and Puri (1988, 1991, 1994) on rank and signed rank methods for autoregressive-moving average (ARMA) and a few other [see, e.g., Hallin and Werker (2003)] time series models. But, except for a componentwise rank approach [Hallin, Ingenbleek and Puri (1989), Hallin and Puri (1995)] to the problem of testing multivariate white noise against vector autoregressive-moving average (VARMA) dependence—suffering the same lack of affine-invariance as its nonserial counterparts—and an affine-invariant approach to the same problem based on interdirections and pseudo-Mahalanobis ranks [Hallin and Paindaveine (2002b)], the multivariate situation so far remains virtually unexplored from this point of view.

Our objective is to develop a complete, fully operational theory of optimal signed rank tests for linear restrictions on the parameters of multiresponse linear models with VARMA errors and unspecified elliptically symmetric innovation densities. These tests are based on the ranks of *pseudo-Mahalanobis distances* and on normed standardized residuals computed from Tyler's (1987) scatter matrix. They generalize the univariate signed rank procedures proposed by Hallin and Puri (1991, 1994) (Tyler-normed residuals playing the role of multivariate signs). Two types of optimality properties are considered, both in the local and asymptotic sense, a la Le Cam: (a) (fixed-score procedures) local asymptotic minimaxity at selected radial densities, and (b) (estimated-score procedures) local asymptotic minimaxity, uniform over the set of all radial densities (satisfying adequate regularity assumptions).

Fulfilling such an objective requires a series of steps, each of which plays an essential role in the construction of the final methods:

(i) defining the adequate (asymptotically sufficient, in the Le Cam sense) rank-based measures of serial dependence, establishing the required asymptotic representation and central-limit results and constructing the optimal tests for fully specified values of the parameters; this is the aim of the present paper, which also works out the algebraic problems related with the singularity of information matrices;

(ii) characterizing the class of linear hypotheses that are invariant under linear transformations, and for which affine-invariant multivariate rank tests make sense (a problem that does not appear in one-dimensional setting); this characterization is obtained in Hallin and Paindaveine (2003);

(iii) establishing the asymptotic linearity of the test statistics we are obtaining in this paper [see (i) above]; this linearity is required if estimated residuals are to be substituted for the exact ones in the computation of multivariate ranks and signs, and is the subject of Hallin and Paindaveine (2004a);

(iv) finally, obtaining the optimal aligned sign and rank tests for linear restrictions, with a detailed and explicit description of some important special cases such as testing for the orders of a VARMA error or a rank-based solution to the multivariate Durbin–Watson problem; see Hallin and Paindaveine (2005).

1.2. *The benefits of a rank-based approach.* Introducing ranks in multivariate time series problems is not just a mathematically challenging exercise, or a matter of theoretical aesthetics. The benefits of rank-based meth-ods indeed are multiple and "mutually orthogonal" in the sense that none of them is obtained at the expense of the others.

If efficiency is the main objective, the generalized Chernoff–Savage and Hodges–Lehmann results of Section 7.2 are important selling points: the fact that asymptotic relative efficiencies with respect to daily-practice methods are never less than one, for instance, is not a small advantage.

But there is much more. The tools we are using here were first developed (in the simpler context of independent observations) in the robustness literature. Robustness (in the vague but reasonably well-understood sense of "resistance to outliers") and efficiency objectives, to a large extent, thus can be attained simultaneously. Invariance—whether strict (with respect to affine transformations) or approximate (with respect to order-preserving radial transformations: see Section 4 for precise statements)—is another fundamental feature of the methods we are proposing. A major consequence of invariance indeed is distribution-freeness, hence exact, similar, unbiased tests, even in the presence of heavy-tailed distributions. The methods we are developing thus remain valid under a very broad class of densities, whereas everyday practice requires finite second-order moments, often fourth-order ones.

Due to the need for consistent estimation of nuisance parameters, some of these benefits of invariance (such as heavy-tail validity) unavoidably have to be tuned down when testing for linear restrictions on the parameters. One way around this problem would consist in modeling median or quantile (auto)regressions, but this is much beyond the scope of the present paper. Most of the nice consequences of invariance however remain, under approximate or asymptotic form.

Finally, the methods proposed are fully applicable; see Hallin and Paindaveine (2004b) for an application to VAR order identification.

1.3. *Outline of the paper.* The starting point in this paper is a local asymptotic normality (LAN) result by Garel and Hallin (1995). This LAN result allows for deriving testing procedures that are locally and asymptotically optimal under a given innovation density $f$, based on a non-Gaussian form of cross-covariances, the *residual $f$-cross-covariance matrices*.

However, due to the possibility of singular local information matrices [such singularities occur as soon as the VARMA$(p_1, q_1)$ neighborhood of a null VARMA$(p_0, q_0)$ model with $p_0 < p_1$ and $q_0 < q_1$ is considered], the optimal test statistics involve unpleasant generalized inverses, which darkens their asymptotic behavior. Therefore, we first restate the LAN property by rewriting the central sequence in a way that explicitly involves the ranks (in the algebraic sense) of local information matrices, and allows for "generalized inverse-free" locally and asymptotically optimal procedures (see the comments after Proposition 1). Next, we show how to replace the "parametric" residual $f$-cross-covariance matrices appearing in the central sequences with "nonparametric" versions, based on the ranks of the Mahalanobis distances and the estimated standardized residuals computed from Tyler's scatter matrix.

Tyler's scatter matrix enjoys highly desirable equivariance/invariance properties. These properties extend to our test statistics; in particular, they are asymptotically invariant under monotone radial transformations of the residuals, hence, asymptotically distribution-free with respect to the underlying radial density. They also are asymptotically distribution-free with respect to the scatter parameter; it should be stressed, however, that this latter property does not follow from any affine-invariance property. Unlike the null hypotheses of location or randomness considered in Hallin and Paindaveine (2002a, b), hypotheses involving general VARMA dependence, as a rule, are not affine-invariant: see Hallin and Paindaveine (2003) for a precise characterization. Actually, our test statistics are strictly affine-invariant whenever the underlying testing problems are, which is of course the most sensible affine-invariance property we can hope for.

We conclude the paper by computing the asymptotic relative efficiencies of the proposed nonparametric procedures with respect to the Gaussian ones, showing that the AREs, as well as the generalized Chernoff–Savage and Hodges–Lehmann theorems obtained in Hallin and Paindaveine (2002b), are still valid here.

The paper is organized as follows. In Section 2 we describe the testing problem under study, and state the main assumptions to be made. The LAN structure of the model is established in Section 3, with a central sequence that exploits the assumptions of elliptical symmetry. The multivariate counterparts of traditional ranks and signs are based on Tyler's scatter matrix, the corresponding *Tyler residuals* and the so-called *pseudo-Mahalanobis ranks*. These concepts are defined in Section 4, where their consistency and invariance properties are also derived. They are used, in Section 5, in the definition of a concept of nonparametric residual cross-covariance matrices, extending to the multivariate context the notion of rank-based autocorrelations developed in Hallin and Puri (1988, 1991). These matrices allow for a reconstruction of central sequences, hence, for nonparametric locally asymptotically optimal procedures. Two types of optimality properties are considered in Section 6, both in the local and asymptotic sense, a la Le Cam (we use the term "minimaxity" even though the tests are "maximin" rather than "minimax"):

(a) (fixed-score procedures) local asymptotic minimaxity at selected radial densities, and

(b) (estimated-score procedures) local asymptotic minimaxity, uniform over a broad class of radial densities.

In both cases the proposed tests remain valid under arbitrary elliptically symmetric innovation densities, including those with infinite variance. In Section 7 the asymptotic relative efficiencies of the proposed procedures, with respect to their Gaussian counterparts (based on classical cross-covariances), are derived. Proofs are concentrated in the Appendix.

## 2. Notation and main assumptions.

Consider the VARMA$(p_1, q_1)$ model defined by the stochastic difference equation

$$\left(\mathbf{I}_k - \sum_{i=1}^{p_1} \mathbf{A}_i L^i\right)\mathbf{X}_t = \left(\mathbf{I}_k + \sum_{i=1}^{q_1} \mathbf{B}_i L^i\right)\boldsymbol{\varepsilon}_t, \qquad t \in \mathbb{Z}, \tag{1}$$

where $\mathbf{A}_1, \ldots, \mathbf{A}_{p_1}, \mathbf{B}_1, \ldots, \mathbf{B}_{q_1}$ are $k \times k$ real matrices ($\mathbf{I}_k$ stands for the $k$-dimensional identity matrix), $L$ denotes the lag operator and $\{\boldsymbol{\varepsilon}_t | t \in \mathbb{Z}\}$ is an absolutely continuous $k$-variate white noise process. The parameter of interest is $\boldsymbol{\theta} := ((\operatorname{vec}\mathbf{A}_1)', \ldots, (\operatorname{vec}\mathbf{A}_{p_1})', (\operatorname{vec}\mathbf{B}_1)', \ldots, (\operatorname{vec}\mathbf{B}_{q_1})')'$, with values in $\mathbb{R}^{k^2(p_1+q_1)}$.

Fixing some value

$$\begin{aligned}\boldsymbol{\theta}_0 := ((\operatorname{vec}\mathbf{A}_1)', \ldots, (\operatorname{vec}\mathbf{A}_{p_0})', \mathbf{0}'_{k^2(p_1-p_0)\times 1}, \\ (\operatorname{vec}\mathbf{B}_1)', \ldots, (\operatorname{vec}\mathbf{B}_{q_0})', \mathbf{0}'_{k^2(q_1-q_0)\times 1})'\end{aligned}$$

of the parameter $\boldsymbol{\theta}$ that satisfies Assumption (A) below, we want to test the null hypothesis $\boldsymbol{\theta} = \boldsymbol{\theta}_0$ against the alternative $\boldsymbol{\theta} \neq \boldsymbol{\theta}_0$. Choosing $p_0 < p_1$ and/or $q_0 < q_1$ allows one to test the adequacy of the specified VARMA coefficients in $\boldsymbol{\theta}_0$, while contemplating the possibility of possibly higher-order VARMA models. If the order is not an issue, one can just let $p_0 = p_1$ and $q_0 = q_1$.

The null VARMA model must satisfy the usual causality and invertibility conditions. More precisely, we assume the following on $\boldsymbol{\theta}_0$:

ASSUMPTION (A). All solutions of $\det(\mathbf{I}_k - \sum_{i=1}^{p_0} \mathbf{A}_i z^i) = 0$, $z \in \mathbb{C}$, and $\det(\mathbf{I}_k + \sum_{i=1}^{q_0} \mathbf{B}_i z^i) = 0$, $z \in \mathbb{C}$ ($|\mathbf{A}_{p_0}| \neq 0 \neq |\mathbf{B}_{q_0}|$), lie outside the unit ball in $\mathbb{C}$. Moreover, the greatest common left divisor of $\mathbf{I}_k - \sum_{i=1}^{p_0} \mathbf{A}_i z^i$ and $\mathbf{I}_k + \sum_{i=1}^{q_0} \mathbf{B}_i z^i$ is $\mathbf{I}_k$.

Write $\mathbf{A}(L)$ and $\mathbf{B}(L)$ for the difference operators $\mathbf{I}_k - \sum_{i=1}^{p_0} \mathbf{A}_i L^i$ and $\mathbf{I}_k + \sum_{i=1}^{q_0} \mathbf{B}_i L^i$, respectively. Letting $\mathbf{B}_0 := \mathbf{I}_k$, recall that the Green's matrices $\mathbf{H}_u$, $u \in \mathbb{N}$, associated with $\mathbf{B}(L)$ are defined by the linear recursion $\sum_{i=0}^{\min(q_0,u)} \mathbf{B}_i \times \mathbf{H}_{u-i} = \delta_{u0}\mathbf{I}_k$, where $\delta_{u0} = 1$ if $u = 0$, and $\delta_{u0} = 0$ otherwise. Assumption (A) also allows for defining these Green's matrices by

$$\sum_{u=0}^{+\infty} \mathbf{H}_u z^u := \left(\mathbf{I}_k + \sum_{i=1}^{q_0} \mathbf{B}_i z^i\right)^{-1}, \qquad z \in \mathbb{C}, |z| < 1; \tag{2}$$

as a consequence, the same matrices also satisfy $\sum_{i=0}^{\min(q_0,u)} \mathbf{B}_i' \mathbf{H}_{u-i}' = \delta_{u0}\mathbf{I}_k$. Assumption (A), moreover, ensures the existence of some $\varepsilon > 0$ such that the series in (2) is absolutely and uniformly convergent in the closed ball in $\mathbb{C}$ with center 0 and radius $1+\varepsilon$. Consequently, $\|\mathbf{H}_u\|(1+\varepsilon)^u$ goes to 0 as $u$ goes to infinity. This exponential decrease ensures that $(\|\mathbf{H}_u\|, u \in \mathbb{N})$ belongs to $l^p(\mathbb{N})$ for all $p > 0$, where $l^p(\mathbb{N})$ denotes the set of all sequences $(x_u, u \in \mathbb{N})$ for which $\sum_{u=0}^{\infty} |x_u|^p < \infty$. Of course, the same remarks also hold, with obvious changes, for Green's matrices $\mathbf{G}_u$, $u \in \mathbb{N}$, associated with the operator $\mathbf{A}(L)$. For simplicity, we do not indicate the strong dependence on $\boldsymbol{\theta}_0$ of $\mathbf{G}_u$ and $\mathbf{H}_u$, which of course are associated with the null operators $\mathbf{A}(L)$ and $\mathbf{B}(L)$.

Under Assumption (A), the white noise $\{\boldsymbol{\varepsilon}_t\}$ is $\{\mathbf{X}_t\}$'s innovation process. The set of Assumptions (B) deals with the density of this innovation. As indicated, we restrict ourselves to a class of elliptically symmetric densities.

ASSUMPTION (B1). Denote by $\boldsymbol{\Sigma}$ a symmetric positive definite $k \times k$ matrix, and by $f : \mathbb{R}_0^+ \to \mathbb{R}^+$ a nonnegative function, such that $f > 0$ a.e. and $\int_0^\infty r^{k-1} f(r)\,dr < \infty$. We will assume throughout that $\{\boldsymbol{\varepsilon}_1^{(n)}, \ldots, \boldsymbol{\varepsilon}_n^{(n)}\}$ is a finite realization of an elliptic white noise process with scatter matrix $\boldsymbol{\Sigma}$ and radial density $f$, that is, such that its probability density at $(\mathbf{z}_1, \ldots, \mathbf{z}_n) \in \mathbb{R}^{nk}$ is of the form $\prod_{t=1}^n \underline{f}(\mathbf{z}_t^{(n)}; \boldsymbol{\Sigma}, f)$, where

$$\underline{f}(\mathbf{z}_1; \boldsymbol{\Sigma}, f) := c_{k,f} (\det \boldsymbol{\Sigma})^{-1/2} f(\|\mathbf{z}_1\|_{\boldsymbol{\Sigma}}), \qquad \mathbf{z}_1 \in \mathbb{R}^k. \tag{3}$$

Here $\|\mathbf{z}\|_{\boldsymbol{\Sigma}} := (\mathbf{z}'\boldsymbol{\Sigma}^{-1}\mathbf{z})^{1/2}$ denotes the norm of $\mathbf{x}$ in the metric associated with $\boldsymbol{\Sigma}$. The constant $c_{k,f}$ is the normalization factor $(\omega_k \mu_{k-1;f})^{-1}$, where $\omega_k$ stands for the $(k-1)$-dimensional Lebesgue measure of the unit sphere $\mathcal{S}^{k-1} \subset \mathbb{R}^k$, and $\mu_{l;f} := \int_0^\infty r^l f(r)\,dr$.

Note that the scatter matrix $\boldsymbol{\Sigma}$ in (3) need not be (a multiple of) the covariance matrix of the observations, which may not exist, and that $f$ is not, strictly speaking, a probability density; see Hallin and Paindaveine (2002a) for a discussion. Moreover, $\boldsymbol{\Sigma}$ and $f$ are identified up to an arbitrary scale transformation only. More precisely, for any $a > 0$, letting $\boldsymbol{\Sigma}_a := a^2\boldsymbol{\Sigma}$ and $f_a(r) := f(ar)$, we have $\underline{f}(\mathbf{x}; \boldsymbol{\Sigma}, f) = \underline{f}(\mathbf{x}; \boldsymbol{\Sigma}_a, f_a)$. This will be harmless in the sequel since cross-covariances, central sequences, as well as all statistics considered, are insensitive to a variation of $a$.

We will denote by $\boldsymbol{\Sigma}^{-1/2}$ the unique upper-triangular $k \times k$ array with positive diagonal elements that satisfies $\boldsymbol{\Sigma}^{-1} = (\boldsymbol{\Sigma}^{-1/2})'\boldsymbol{\Sigma}^{-1/2}$. With this notation, $\boldsymbol{\Sigma}^{-1/2}\boldsymbol{\varepsilon}_1^{(n)}/\|\boldsymbol{\Sigma}^{-1/2}\boldsymbol{\varepsilon}_1^{(n)}\|, \ldots, \boldsymbol{\Sigma}^{-1/2}\boldsymbol{\varepsilon}_n^{(n)}/\|\boldsymbol{\Sigma}^{-1/2}\boldsymbol{\varepsilon}_n^{(n)}\|$ are i.i.d., and uniformly distributed over $\mathcal{S}^{k-1}$. Similarly, $\|\boldsymbol{\Sigma}^{-1/2}\boldsymbol{\varepsilon}_1^{(n)}\|, \ldots, \|\boldsymbol{\Sigma}^{-1/2}\boldsymbol{\varepsilon}_n^{(n)}\|$ are i.i.d. with probability density function

$$\tilde{f}_k(r) := (\mu_{k-1;f})^{-1} r^{k-1} f(r) I_{[r>0]}, \tag{4}$$

where $I_E$ denotes the indicator function associated with the Borel set $E$. The terminology *radial density* will be used for $f$ and $\tilde{f}_k$ indifferently (though only $\tilde{f}_k$, of course, is a genuine probability density). We denote by $\tilde{F}_k$ the distribution function associated with $\tilde{f}_k$.

Write $\mathcal{H}^{(n)}(\boldsymbol{\theta}_0, \boldsymbol{\Sigma}, f)$ for the hypothesis under which an observation $\mathbf{X}^{(n)} := (\mathbf{X}_1^{(n)}, \ldots, \mathbf{X}_n^{(n)})'$ is generated by the VARMA$(p_0, q_0)$ model (1) with parameter value $\boldsymbol{\theta}_0$ satisfying Assumption (A) and innovation process satisfying Assumption (B1) with "scatter parameter" $\boldsymbol{\Sigma}$ and radial density $f$. Our objective is to test $\mathcal{H}^{(n)}(\boldsymbol{\theta}_0) := \bigcup_{\boldsymbol{\Sigma}} \bigcup_f \mathcal{H}^{(n)}(\boldsymbol{\theta}_0, \boldsymbol{\Sigma}, f)$ against $\bigcup_{\boldsymbol{\theta} \neq \boldsymbol{\theta}_0} \mathcal{H}^{(n)}(\boldsymbol{\theta})$. Consequently, $\boldsymbol{\Sigma}$ and $f$ play the role of nuisance parameters; note that the unions, in the definition of $\mathcal{H}^{(n)}(\boldsymbol{\theta}_0)$, are taken over *all* possible values of $\boldsymbol{\Sigma}$ and $f$.

The methodology we will adopt in the derivation of optimality results is based on Le Cam's asymptotic theory of experiments. This requires the model (characterized by some fixed value of $\boldsymbol{\Sigma}$ and $f$) at which optimality is sought to be locally and asymptotically normal (LAN). LAN, of course, does not hold without a few, rather mild, regularity conditions: finite second-order moments and finite Fisher information for $\underline{f}$, and quadratic mean differentiability of $\underline{f}^{1/2}$. These technical assumptions are taken into account in Assumptions (B1$'$) and (B2), in a form that is adapted to the elliptical context. We insist, however, on the fact that these assumptions are made on the density at which optimality is desired, not on the actual density.

Assumption (B1$'$). Same as Assumption (B1), but we further assume that $\mu_{k+1;f} < \infty$.

Note that Assumption (B1$'$) is also required when Gaussian procedures are to be considered, since these also require the second-order moments of the underlying distribution to be finite.

Assumption (B2) is strictly equivalent to the assumption that $\underline{f}^{1/2}$ is differentiable in quadratic mean [see Hallin and Paindaveine (2002a), Proposition 1]. However, it has the important advantage of involving univariate quadratic mean differentiability only. Let $L^2(\mathbb{R}_0^+, \mu_l)$ denote the space of all functions that are square-integrable with respect to the Lebesgue measure with weight $r^l$ on $\mathbb{R}_0^+$, that is, the space of measurable functions $h : \mathbb{R}_0^+ \to \mathbb{R}$ satisfying $\int_0^\infty [h(r)]^2 r^l \, dr < \infty$.

Assumption (B2). The square root $f^{1/2}$ of the radial density $f$ is in $W^{1,2}(\mathbb{R}_0^+, \mu_{k-1})$, where $W^{1,2}(\mathbb{R}_0^+, \mu_{k-1})$ denotes the subspace of $L^2(\mathbb{R}_0^+, \mu_{k-1})$ containing all functions admitting a weak derivative that also belongs to $L^2(\mathbb{R}_0^+, \mu_{k-1})$.

Denoting by $(f^{1/2})'$ the weak derivative of $f^{1/2}$ in $L^2(\mathbb{R}_0^+, \mu_{k-1})$, let $\varphi_f := -2\frac{(f^{1/2})'}{f^{1/2}}$: Assumption (B2) ensures the finiteness of the *radial Fisher information*

$$\mathcal{I}_{k,f} := (\mu_{k-1;f})^{-1} \int_0^\infty [\varphi_f(r)]^2 r^{k-1} f(r)\, dr.$$

In the last set of assumptions, we describe the score functions $K_1$, $K_2$ to be used when building rank-based statistics in this serial context.

ASSUMPTION (C). The score functions $K_\ell : ]0,1[ \to \mathbb{R}$, $\ell = 1, 2$, are continuous differences of two monotone increasing functions, and satisfy $\int_0^1 [K_\ell(u)]^2\, du < \infty$ ($\ell = 1, 2$).

The score functions yielding locally and asymptotically optimal procedures, as we shall see, are of the form $K_1 := \varphi_{f_\star} \circ \tilde{F}_{\star k}^{-1}$ and $K_2 := \tilde{F}_{\star k}^{-1}$, for some radial density $f_\star$. Assumption (C) then takes the form of an assumption on $f_\star$:

ASSUMPTION (C$'$). The radial density $f_\star$ is such that $\varphi_{f_\star}$ is the continuous difference of two monotone increasing functions, $\mu_{k+1;f_\star} < \infty$ and $\int_0^\infty [\varphi_{f_\star}(r)]^2 r^{k-1} f_\star(r)\, dr < \infty$.

The assumption being the difference of two monotone functions, which characterizes the functions with bounded variation, is extremely mild. In most cases ($f_\star$ normal, double exponential, ...), $\varphi_{f_\star}$ itself is monotone increasing, and, without loss of generality, this will be assumed to hold for the proofs. The multivariate $t$-distributions, however, provide examples of nonmonotone score functions $\varphi_{f_\star}$ satisfying Assumption (C$'$).

**3. Local asymptotic normality.** Let $\mathbf{A}(L)$ and $\mathbf{B}(L)$ be such that $\mathbf{A}_i := \mathbf{0}$ for $i = p_0 + 1, \ldots, p_1$, and $\mathbf{B}_i := \mathbf{0}$ for $i = q_0 + 1, \ldots, q_1$, and consider the sequences of linear difference operators

$$\mathbf{A}^{(n)}(L) := \mathbf{I}_k - \sum_{i=1}^{p_1} (\mathbf{A}_i + n^{-1/2} \boldsymbol{\gamma}_i^{(n)}) L^i \quad \text{and}$$

$$\mathbf{B}^{(n)}(L) := \mathbf{I}_k + \sum_{i=1}^{q_1} (\mathbf{B}_i + n^{-1/2} \boldsymbol{\delta}_i^{(n)}) L^i,$$

where the vector $\boldsymbol{\tau}^{(n)} := ((\operatorname{vec} \boldsymbol{\gamma}_1^{(n)})', \ldots, (\operatorname{vec} \boldsymbol{\gamma}_{p_1}^{(n)})', (\operatorname{vec} \boldsymbol{\delta}_1^{(n)})', \ldots, (\operatorname{vec} \boldsymbol{\delta}_{q_1}^{(n)})')' \in \mathbb{R}^{k^2(p_1+q_1)}$ is such that $\sup_n (\boldsymbol{\tau}^{(n)})' \boldsymbol{\tau}^{(n)} < \infty$. These operators define a sequence of VARMA models

$$\mathbf{A}^{(n)}(L) \mathbf{X}_t = \mathbf{B}^{(n)}(L) \boldsymbol{\varepsilon}_t, \qquad t \in \mathbb{Z},$$

hence, in the notation of Section 2, the sequence of local alternatives $\mathcal{H}^{(n)}(\boldsymbol{\theta}_0 + n^{-1/2}\boldsymbol{\tau}^{(n)}, \boldsymbol{\Sigma}, f)$.

Let $(\mathbf{Z}_1^{(n)}(\boldsymbol{\theta}_0), \ldots, \mathbf{Z}_n^{(n)}(\boldsymbol{\theta}_0))$ be the $n$-tuple of residuals computed from the initial values $\boldsymbol{\varepsilon}_{-q_0+1}, \ldots, \boldsymbol{\varepsilon}_0$ and $\mathbf{X}_{-p_0+1}^{(n)}, \ldots, \mathbf{X}_0^{(n)}$ and the observed series $(\mathbf{X}_1^{(n)}, \ldots, \mathbf{X}_n^{(n)})$ via the relation

$$\mathbf{Z}_t^{(n)}(\boldsymbol{\theta}_0) = \sum_{i=0}^{t-1}\sum_{j=0}^{p_0} \mathbf{H}_i \mathbf{A}_j \mathbf{X}_{t-i-j}^{(n)} + (\mathbf{H}_{t+q_0-1} \cdots \mathbf{H}_t) \begin{pmatrix} \mathbf{I}_k & \mathbf{0} & \ldots & \mathbf{0} \\ \mathbf{B}_1 & \mathbf{I}_k & \ldots & \mathbf{0} \\ \vdots & \vdots & \ddots & \vdots \\ \mathbf{B}_{q_0-1} & \mathbf{B}_{q_0-2} & \ldots & \mathbf{I}_k \end{pmatrix} \begin{pmatrix} \boldsymbol{\varepsilon}_{-q_0+1} \\ \vdots \\ \boldsymbol{\varepsilon}_0 \end{pmatrix}.$$

Assumption (A) ensures that neither the (generally unobserved) values $(\boldsymbol{\varepsilon}_{-q_0+1}, \ldots, \boldsymbol{\varepsilon}_0)$ of the innovation, nor the initial values $(\mathbf{X}_{-p_0+1}^{(n)}, \ldots, \mathbf{X}_0^{(n)})$, have an influence on asymptotic results, so that they all safely can be set to zero in the sequel. Decompose $\mathbf{Z}_t^{(n)}(\boldsymbol{\theta}_0)$ into

$$\mathbf{Z}_t^{(n)}(\boldsymbol{\theta}_0) = d_t^{(n)}(\boldsymbol{\theta}_0, \boldsymbol{\Sigma})\, \boldsymbol{\Sigma}^{1/2} \mathbf{U}_t^{(n)}(\boldsymbol{\theta}_0, \boldsymbol{\Sigma}),$$

where $d_t^{(n)}(\boldsymbol{\theta}_0, \boldsymbol{\Sigma}) := \|\mathbf{Z}_t^{(n)}(\boldsymbol{\theta}_0)\|_{\boldsymbol{\Sigma}}$ and $\mathbf{U}_t^{(n)}(\boldsymbol{\theta}_0, \boldsymbol{\Sigma}) := \boldsymbol{\Sigma}^{-1/2}\mathbf{Z}_t^{(n)}(\boldsymbol{\theta}_0)/d_t^{(n)}(\boldsymbol{\theta}_0, \boldsymbol{\Sigma})$. Writing $\varphi_{\underline{f}} := -2(\mathbf{D}\underline{f}^{1/2})/\underline{f}^{1/2}$, where $\mathbf{D}\underline{f}^{1/2}$ denotes the quadratic mean gradient of $\underline{f}^{1/2}$, define, as in Garel and Hallin (1995), the residual $\underline{f}$-*cross covariance matrix* of lag $i$ as

$$\boldsymbol{\Gamma}_{i;\underline{f}}^{(n)}(\boldsymbol{\theta}_0) := (n-i)^{-1} \sum_{t=i+1}^{n} \varphi_{\underline{f}}(\mathbf{Z}_t^{(n)}(\boldsymbol{\theta}_0))\mathbf{Z}_{t-i}^{(n)\prime}(\boldsymbol{\theta}_0). \tag{5}$$

Due to the elliptical structure of $\underline{f}$, these cross-covariance matrices take the form

$$\boldsymbol{\Gamma}_{i;\boldsymbol{\Sigma},f}^{(n)}(\boldsymbol{\theta}_0) = \frac{1}{n-i}\boldsymbol{\Sigma}'^{-1/2}\Bigg(\sum_{t=i+1}^{n} \varphi_f(d_t^{(n)}(\boldsymbol{\theta}_0, \boldsymbol{\Sigma}))\, d_{t-i}^{(n)}(\boldsymbol{\theta}_0, \boldsymbol{\Sigma}) \times \mathbf{U}_t^{(n)}(\boldsymbol{\theta}_0, \boldsymbol{\Sigma})\mathbf{U}_{t-i}^{(n)\prime}(\boldsymbol{\theta}_0, \boldsymbol{\Sigma})\Bigg)\boldsymbol{\Sigma}'^{1/2}. \tag{6}$$

Hallin and Paindaveine (2002b) are developing optimal procedures based on nonparametric versions of the cross-covariances (5) or (6) for the problem of testing elliptic white noise against VARMA dependence.

Garel and Hallin (1995) established LAN in this setting (in fact, in a more general, possibly nonelliptical, context). The quadratic form in their second-order approximation of local log-likelihoods (hence, also in ours) is not of full rank, due to the well-known singularity of the information matrix for ARMA models. This singularity, quite understandably, has to be taken into account in the derivation of locally optimal inference procedures; in hypothesis testing, this can be achieved in a straightforward way provided the information matrix is factorized in an adequate way (see the comments after Proposition 1 for details). This factorization is not provided by Garel and Hallin (1995), since they just deal with LAN, and not with optimal inference, but it is needed here. Therefore, in Proposition 1, we first restate LAN under a slightly different form, in the spirit of the univariate results of Hallin and Puri (1994). As usual, though, the multivariate case is a bit more intricate, and requires some notational and algebraic preparation.

Associated with any $k$-dimensional linear difference operator of the form $\mathbf{C}(L) := \sum_{i=0}^{\infty} \mathbf{C}_i L^i$ (letting $\mathbf{C}_i = \mathbf{0}$ for $i > s$, this includes, of course, operators of finite order $s$), define for any integers $u$ and $v$ the $k^2u \times k^2v$ matrices

$$\mathbf{C}_{u,v}^{(l)} := \begin{pmatrix} \mathbf{C}_0 \otimes \mathbf{I}_k & \mathbf{0} & \dots & \mathbf{0} \\ \mathbf{C}_1 \otimes \mathbf{I}_k & \mathbf{C}_0 \otimes \mathbf{I}_k & \dots & \mathbf{0} \\ \vdots & & \ddots & \vdots \\ \mathbf{C}_{v-1} \otimes \mathbf{I}_k & \mathbf{C}_{v-2} \otimes \mathbf{I}_k & \dots & \mathbf{C}_0 \otimes \mathbf{I}_k \\ \vdots & & & \vdots \\ \mathbf{C}_{u-1} \otimes \mathbf{I}_k & \mathbf{C}_{u-2} \otimes \mathbf{I}_k & \dots & \mathbf{C}_{u-v} \otimes \mathbf{I}_k \end{pmatrix}$$

and

$$\mathbf{C}_{u,v}^{(r)} := \begin{pmatrix} \mathbf{I}_k \otimes \mathbf{C}_0 & \mathbf{0} & \dots & \mathbf{0} \\ \mathbf{I}_k \otimes \mathbf{C}_1 & \mathbf{I}_k \otimes \mathbf{C}_0 & \dots & \mathbf{0} \\ \vdots & & \ddots & \vdots \\ \mathbf{I}_k \otimes \mathbf{C}_{v-1} & \mathbf{I}_k \otimes \mathbf{C}_{v-2} & \dots & \mathbf{I}_k \otimes \mathbf{C}_0 \\ \vdots & & & \vdots \\ \mathbf{I}_k \otimes \mathbf{C}_{u-1} & \mathbf{I}_k \otimes \mathbf{C}_{u-2} & \dots & \mathbf{I}_k \otimes \mathbf{C}_{u-v} \end{pmatrix},$$

respectively; write $\mathbf{C}_u^{(l)}$ for $\mathbf{C}_{u,u}^{(l)}$ and $\mathbf{C}_u^{(r)}$ for $\mathbf{C}_{u,u}^{(r)}$. With this notation, note that $\mathbf{G}_u^{(l)}$, $\mathbf{G}_u^{(r)}$, $\mathbf{H}_u^{(l)}$ and $\mathbf{H}_u^{(r)}$ are the inverses of $\mathbf{A}_u^{(l)}$, $\mathbf{A}_u^{(r)}$, $\mathbf{B}_u^{(l)}$ and $\mathbf{B}_u^{(r)}$, respectively. Denoting by $\mathbf{C}_{u,v}^{\prime(l)}$ and $\mathbf{C}_{u,v}^{\prime(r)}$ the matrices associated with the transposed operator $\mathbf{C}'(L) := \sum_{i=0}^{\infty} \mathbf{C}'_i L^i$, we also have $\mathbf{G}_u^{\prime(l)} = (\mathbf{A}_u^{\prime(l)})^{-1}$, $\mathbf{H}_u^{\prime(l)} = (\mathbf{B}_u^{\prime(l)})^{-1}$, and so on.

Let $\pi := \max(p_1 - p_0, q_1 - q_0)$ and $\pi_0 := \pi + p_0 + q_0$, and define the $k^2\pi_0 \times k^2(p_1 + q_1)$ matrix $\mathbf{M}_{\boldsymbol{\theta}_0} := (\mathbf{G}_{\pi_0,p_1}^{\prime(l)} \vdots \mathbf{H}_{\pi_0,q_1}^{\prime(l)})$: under Assumption (A) $\mathbf{M}_{\boldsymbol{\theta}_0}$ is of full rank.

Finally, consider the operator $\mathbf{D}(L) := \mathbf{I}_k + \sum_{i=1}^{p_0+q_0} \mathbf{D}_i L^i$ [just as $\mathbf{M}_{\boldsymbol{\theta}_0}$, $\mathbf{D}(L)$ and most quantities defined below depend on $\boldsymbol{\theta}_0$; for simplicity, however, we are dropping this reference to $\boldsymbol{\theta}_0$], where, putting $\mathbf{G}_{-1} = \mathbf{G}_{-2} =$

$\cdots = \mathbf{G}_{-p_0+1} = \mathbf{0} = \mathbf{H}_{-1} = \mathbf{H}_{-2} = \cdots = \mathbf{H}_{-q_0+1}$,

$$\begin{pmatrix} \mathbf{D}'_1 \\ \vdots \\ \mathbf{D}'_{p_0+q_0} \end{pmatrix} := - \begin{pmatrix} \mathbf{G}_{q_0} & \mathbf{G}_{q_0-1} & \cdots & \mathbf{G}_{-p_0+1} \\ \mathbf{G}_{q_0+1} & \mathbf{G}_{q_0} & \cdots & \mathbf{G}_{-p_0+2} \\ \vdots & & \ddots & \vdots \\ \mathbf{G}_{p_0+q_0-1} & \mathbf{G}_{p_0+q_0-2} & \cdots & \mathbf{G}_0 \\ \mathbf{H}_{p_0} & \mathbf{H}_{p_0-1} & \cdots & \mathbf{H}_{-q_0+1} \\ \mathbf{H}_{p_0+1} & \mathbf{H}_{p_0} & \cdots & \mathbf{H}_{-q_0+2} \\ \vdots & & \ddots & \vdots \\ \mathbf{H}_{p_0+q_0-1} & \mathbf{H}_{p_0+q_0-2} & \cdots & \mathbf{H}_0 \end{pmatrix}^{-1} \begin{pmatrix} \mathbf{G}_{q_0+1} \\ \vdots \\ \mathbf{G}_{p_0+q_0} \\ \mathbf{H}_{p_0+1} \\ \vdots \\ \mathbf{H}_{p_0+q_0} \end{pmatrix}.$$

Note that $\mathbf{D}(L)\mathbf{G}'_t = \mathbf{0}$ for $t = q_0 + 1, \ldots, p_0 + q_0$, and $\mathbf{D}(L)\mathbf{H}'_t = \mathbf{0}$ for $t = p_0 + 1, \ldots, p_0 + q_0$.

Let $\{\mathbf{\Psi}_t^{(1)}, \ldots, \mathbf{\Psi}_t^{(p_0+q_0)}\}$ be a set of $k \times k$ matrices forming a fundamental system of solutions of the homogeneous linear difference equation associated with $\mathbf{D}(L)$. Such a system can be obtained, for instance, from Green's matrices of the operator $\mathbf{D}(L)$ [see, e.g., Hallin (1986)]. Defining

$$\bar{\mathbf{\Psi}}_m(\boldsymbol{\theta}_0) := \begin{pmatrix} \mathbf{\Psi}_{\pi+1}^{(1)} & \cdots & \mathbf{\Psi}_{\pi+1}^{(p_0+q_0)} \\ \mathbf{\Psi}_{\pi+2}^{(1)} & \cdots & \mathbf{\Psi}_{\pi+2}^{(p_0+q_0)} \\ \vdots & & \vdots \\ \mathbf{\Psi}_m^{(1)} & \cdots & \mathbf{\Psi}_m^{(p_0+q_0)} \end{pmatrix} \otimes \mathbf{I}_k \qquad (m > \pi),$$

the Casorati matrix $\mathbf{C}_\Psi$ associated with $\mathbf{D}(L)$ is $\bar{\mathbf{\Psi}}_{\pi_0}$. Putting

$$\mathbf{P}_{\boldsymbol{\theta}_0} := \begin{pmatrix} \mathbf{I}_{k^2\pi} & \mathbf{0} \\ \mathbf{0} & \mathbf{C}_\Psi^{-1} \end{pmatrix} \quad \text{and} \quad \mathbf{Q}_{\boldsymbol{\theta}_0}^{(n)} := \mathbf{H}_{n-1}^{(r)} \mathbf{B}_{n-1}^{\prime(l)} \begin{pmatrix} \mathbf{I}_{k^2\pi} & \mathbf{0} \\ \mathbf{0} & \bar{\mathbf{\Psi}}_{n-1} \end{pmatrix},$$

let

$$\begin{aligned} \mathbf{S}_{\mathbf{\Sigma},f}^{(n)}(\boldsymbol{\theta}_0) &:= ((n-1)^{1/2}(\operatorname{vec} \mathbf{\Gamma}_{1;\mathbf{\Sigma},f}^{(n)}(\boldsymbol{\theta}_0))', \ldots, \\ &\qquad (n-i)^{1/2}(\operatorname{vec} \mathbf{\Gamma}_{i;\mathbf{\Sigma},f}^{(n)}(\boldsymbol{\theta}_0))', \ldots, (\operatorname{vec} \mathbf{\Gamma}_{n-1;\mathbf{\Sigma},f}^{(n)}(\boldsymbol{\theta}_0))')', \\ n^{1/2}\mathbf{T}_{\mathbf{\Sigma},f}^{(n)}(\boldsymbol{\theta}_0) &:= \mathbf{Q}_{\boldsymbol{\theta}_0}^{(n)\prime} \mathbf{S}_{\mathbf{\Sigma},f}^{(n)}(\boldsymbol{\theta}_0) \end{aligned}$$

and

$$\mathbf{J}_{\boldsymbol{\theta}_0,\mathbf{\Sigma}} := \lim_{n \to +\infty} \mathbf{Q}_{\boldsymbol{\theta}_0}^{(n)\prime} [\mathbf{I}_{n-1} \otimes (\mathbf{\Sigma} \otimes \mathbf{\Sigma}^{-1})] \mathbf{Q}_{\boldsymbol{\theta}_0}^{(n)} \tag{7}$$

[convergence in (7) follows from the exponential decrease, as $u \to \infty$, of Green's matrices $\mathbf{G}_u$ and $\mathbf{H}_u$; see the comment after Assumption (A)]. Local asymptotic normality, for fixed $\mathbf{\Sigma}$ and $f$, of the model described in Section 2 then can be stated in the following way.

PROPOSITION 1. *Assume that Assumptions* (A), (B1′) *and* (B2) *hold. Then the logarithm* $L^{(n)}_{\boldsymbol{\theta}_0+n^{-1/2}\boldsymbol{\tau}^{(n)}/\boldsymbol{\theta}_0;\boldsymbol{\Sigma},f}$ *of the likelihood ratio of* $\mathcal{H}^{(n)}(\boldsymbol{\theta}_0+n^{-1/2}\boldsymbol{\tau}^{(n)},\boldsymbol{\Sigma},f)$ *with respect to* $\mathcal{H}^{(n)}(\boldsymbol{\theta}_0,\boldsymbol{\Sigma},f)$ *is such that*

$$L^{(n)}_{\boldsymbol{\theta}_0+n^{-1/2}\boldsymbol{\tau}^{(n)}/\boldsymbol{\theta}_0;\boldsymbol{\Sigma},f}(\mathbf{X}^{(n)})=(\boldsymbol{\tau}^{(n)})'\boldsymbol{\Delta}^{(n)}_{\boldsymbol{\Sigma},f}(\boldsymbol{\theta}_0)-\tfrac{1}{2}(\boldsymbol{\tau}^{(n)})'\boldsymbol{\Gamma}_{\boldsymbol{\Sigma},f}(\boldsymbol{\theta}_0)\boldsymbol{\tau}^{(n)}+o_{\mathrm{P}}(1),$$

*as* $n\to\infty$, *under* $\mathcal{H}^{(n)}(\boldsymbol{\theta}_0,\boldsymbol{\Sigma},f)$, *with the* central sequence

$$\boldsymbol{\Delta}^{(n)}_{\boldsymbol{\Sigma},f}(\boldsymbol{\theta}_0):=n^{1/2}\mathbf{M}'_{\boldsymbol{\theta}_0}\mathbf{P}'_{\boldsymbol{\theta}_0}\mathbf{T}^{(n)}_{\boldsymbol{\Sigma},f}(\boldsymbol{\theta}_0),\tag{8}$$

*and with the asymptotic* information matrix $\boldsymbol{\Gamma}_{\boldsymbol{\Sigma},f}(\boldsymbol{\theta}_0):=\frac{\mu_{k+1;f}\mathcal{I}_{k,f}}{k^2\mu_{k-1;f}}\mathbf{N}_{\boldsymbol{\theta},\boldsymbol{\Sigma}}$, *where*

$$\mathbf{N}_{\boldsymbol{\theta},\boldsymbol{\Sigma}}:=\mathbf{M}'_{\boldsymbol{\theta}_0}\mathbf{P}'_{\boldsymbol{\theta}_0}\mathbf{J}_{\boldsymbol{\theta}_0,\boldsymbol{\Sigma}}\mathbf{P}_{\boldsymbol{\theta}_0}\mathbf{M}_{\boldsymbol{\theta}_0}.\tag{9}$$

*Moreover,* $\boldsymbol{\Delta}^{(n)}_{\boldsymbol{\Sigma},f}(\boldsymbol{\theta}_0)$, *still under* $\mathcal{H}^{(n)}(\boldsymbol{\theta}_0,\boldsymbol{\Sigma},f)$, *is asymptotically* $\mathcal{N}_{k^2(p_1+q_1)}(\mathbf{0},\boldsymbol{\Gamma}_{\boldsymbol{\Sigma},f}(\boldsymbol{\theta}_0))$.

For the proof see the Appendix.

The benefits of expressing $\boldsymbol{\Delta}^{(n)}_{\boldsymbol{\Sigma},f}$ and $\boldsymbol{\Gamma}_{\boldsymbol{\Sigma},f}$ as in (8) and (9) stem from the following elementary facts. Sequences of local experiments under LAN converge, in the Le Cam sense, to Gaussian shift experiments, so that optimal tests for the limit Gaussian shifts determine the form of locally asymptotically optimal tests for the original problem. Consider the problem of testing $\mathcal{H}_0:\boldsymbol{\tau}=\mathbf{0}$ against $\mathcal{H}_1:\boldsymbol{\tau}\neq\mathbf{0}$ in the single-observation $\ell$-variate Gaussian shift experiment under which $\boldsymbol{\Delta}\sim\mathcal{N}_\ell(\boldsymbol{\Gamma\tau},\boldsymbol{\Gamma})$. Let $m:=\operatorname{rank}(\boldsymbol{\Gamma})$. If $m=\ell$, the optimal ($\alpha$-level maximin) test consists in rejecting $\mathcal{H}_0$ for large values of $\boldsymbol{\Delta}'\boldsymbol{\Gamma}^{-1}\boldsymbol{\Delta}$, the null distribution of which is $\chi^2_\ell$. Whenever $m<\ell$, $\boldsymbol{\Gamma}$ is singular, and this does not hold anymore. However, if we succeed in writing $\boldsymbol{\Delta}$ and $\boldsymbol{\Gamma}$ in the form

$$\boldsymbol{\Delta}=\mathbf{M}'\mathring{\boldsymbol{\Delta}}\quad\text{and}\quad\boldsymbol{\Gamma}=\mathbf{M}'\mathring{\boldsymbol{\Gamma}}\mathbf{M},\tag{10}$$

where both the $m\times\ell$ matrix $\mathbf{M}$ and the $m\times m$ matrix $\mathring{\boldsymbol{\Gamma}}$ have full rank, the problem of testing $\mathcal{H}_0:\boldsymbol{\tau}=\mathbf{0}$ in the singular $\ell$-variate Gaussian shift experiment for $\boldsymbol{\Delta}\sim\mathcal{N}_\ell(\boldsymbol{\Gamma\tau},\boldsymbol{\Gamma})$ is strictly the same as that of testing $\mathcal{H}_0:\mathbf{M}\boldsymbol{\tau}=\mathbf{0}$ in the full-rank $m$-variate Gaussian shift experiment under which $\mathring{\boldsymbol{\Delta}}\sim\mathcal{N}_m(\mathring{\boldsymbol{\Gamma}}\mathbf{M}\boldsymbol{\tau},\mathring{\boldsymbol{\Gamma}})$. It follows that the optimal ($\alpha$-level maximin) test for $\mathcal{H}_0$ rejects the null hypothesis for large values of $\mathring{\boldsymbol{\Delta}}'\mathring{\boldsymbol{\Gamma}}^{-1}\mathring{\boldsymbol{\Delta}}$, which is $\chi^2_m$ under $\mathcal{H}_0$. Now, clearly, $\mathring{\boldsymbol{\Delta}}'\mathring{\boldsymbol{\Gamma}}^{-1}\mathring{\boldsymbol{\Delta}}=\boldsymbol{\Delta}'\boldsymbol{\Gamma}^{-}\boldsymbol{\Delta}$, where $\boldsymbol{\Gamma}^{-}$ denotes an arbitrary generalized inverse of $\boldsymbol{\Gamma}$, so that, if we succeed in writing $\boldsymbol{\Delta}$ and $\boldsymbol{\Gamma}$ in the form (10), the somewhat unpleasant recourse to generalized inverses is not required anymore. This is exactly what expressions (8) and (9) are allowing for. As a consequence, the degeneracy of the information matrix is kind of hidden in the explicit forms of the optimal test statistics in Propositions 4, 6 and 7.

## 4. Multivariate ranks and signs: invariance and equivariance.

4.1. *Pseudo-Mahalanobis distances and Tyler residuals.* Likelihoods—hence, the central sequences (8)—are measurable, jointly, with respect to two types of statistics:

(i) the distances $d_t^{(n)}(\boldsymbol{\theta}_0, \boldsymbol{\Sigma})$ between standardized residuals $\boldsymbol{\Sigma}^{-1/2}\mathbf{Z}_t^{(n)}(\boldsymbol{\theta}_0)$ and the origin in $\mathbb{R}^k$, and

(ii) the *normalized standardized residuals* $\mathbf{U}_t^{(n)}(\boldsymbol{\theta}_0, \boldsymbol{\Sigma}) := \boldsymbol{\Sigma}^{-1/2}\mathbf{Z}_t^{(n)}(\boldsymbol{\theta}_0)/d_t^{(n)}(\boldsymbol{\theta}_0, \boldsymbol{\Sigma})$.

The (univariate) distances $d_t^{(n)}(\boldsymbol{\theta}_0, \boldsymbol{\Sigma})$ are i.i.d. over the positive real line, with density (4); their ranks thus have the same distribution-freeness and maximal invariance properties as those of the absolute values of any univariate symmetrically distributed $n$-tuple. The normed standardized residuals $\mathbf{U}_t^{(n)}(\boldsymbol{\theta}_0, \boldsymbol{\Sigma})$ under $\mathcal{H}^{(n)}(\boldsymbol{\theta}_0, \boldsymbol{\Sigma}, f)$ are uniformly distributed over the unit sphere, and, hence, can be viewed as multivariate generalizations of signs.

Unfortunately, both $d_t^{(n)}(\boldsymbol{\theta}_0, \boldsymbol{\Sigma})$ and $\mathbf{U}_t^{(n)}(\boldsymbol{\theta}_0, \boldsymbol{\Sigma})$ involve, in a crucial way, the shape parameter $\boldsymbol{\Sigma}$, which, in practice, is never specified, and has to be estimated from the observations. If the actual underlying distribution has finite second-order moments [i.e., under Assumption (B1′)], a "natural" consistent candidate for estimating $\boldsymbol{\Sigma}$ is the empirical covariance matrix $n^{-1}\sum_{t=1}^{n}\mathbf{Z}_t^{(n)}(\boldsymbol{\theta}_0)(\mathbf{Z}_t^{(n)}(\boldsymbol{\theta}_0))'$. Finite second-order moments, however, are too strong a requirement, as we would like to build testing procedures that are optimal under the assumptions of Proposition 1, but remain valid under much milder conditions, including the case of infinite variances. This rules out the empirical covariance as an estimate of $\boldsymbol{\Sigma}$ and, under the weaker Assumption (B1), which does not require anything about the moments of the underlying distribution, we propose to use Tyler's estimator of scatter [see Tyler (1987)].

This estimator is defined as follows. For any $n$-tuple $\mathbf{Z}^{(n)} := (\mathbf{Z}_1^{(n)}, \mathbf{Z}_2^{(n)}, \ldots, \mathbf{Z}_n^{(n)})$ of $k$-dimensional vectors, denote by $\mathbf{C}_{\mathrm{Tyl}}^{(n)} := \mathbf{C}_{\mathrm{Tyl}}^{(n)}(\mathbf{Z}^{(n)})$ the [unique for $n > k$ $(k-1)$] upper triangular $k \times k$ matrix with positive diagonal elements and a "1" in the upper left corner that satisfies

$$\frac{1}{n}\sum_{t=1}^{n}\left(\frac{\mathbf{C}_{\mathrm{Tyl}}^{(n)}\mathbf{Z}_t^{(n)}}{\|\mathbf{C}_{\mathrm{Tyl}}^{(n)}\mathbf{Z}_t^{(n)}\|}\right)\left(\frac{\mathbf{C}_{\mathrm{Tyl}}^{(n)}\mathbf{Z}_t^{(n)}}{\|\mathbf{C}_{\mathrm{Tyl}}^{(n)}\mathbf{Z}_t^{(n)}\|}\right)' = \frac{1}{k}\mathbf{I}_k; \tag{11}$$

Tyler's estimator of scatter is defined as $\widehat{\boldsymbol{\Sigma}}^{(n)} := (\mathbf{C}_{\mathrm{Tyl}}^{(n)\prime}\mathbf{C}_{\mathrm{Tyl}}^{(n)})^{-1}$.

When computed from the $n$-tuple of residuals $\mathbf{Z}_t^{(n)}(\boldsymbol{\theta}_0)$, $t = 1, \ldots, n$, Tyler's estimator is root-$n$ consistent, up to a multiplicative factor, for the shape matrix $\boldsymbol{\Sigma}$. More precisely, there exists a positive real $a$ such that $\sqrt{n}(\widehat{\boldsymbol{\Sigma}}^{(n)} - a\boldsymbol{\Sigma})$ is $O_{\mathrm{P}}(1)$ as $n \to \infty$ under $\bigcup_f \mathcal{H}^{(n)}(\boldsymbol{\theta}_0, \boldsymbol{\Sigma}, f)$. Tyler's estimator is clearly invariant under permutations of the residuals $\mathbf{Z}_t^{(n)}(\boldsymbol{\theta}_0)$. Moreover, $\mathbf{C}_{\mathrm{Tyl}}^{(n)}$ is strictly affine-equivariant, since

$$\mathbf{C}_{\mathrm{Tyl}}^{(n)}(\mathbf{M}\mathbf{Z}^{(n)}) = d\mathbf{O}\mathbf{C}_{\mathrm{Tyl}}^{(n)}\mathbf{M}^{-1} \tag{12}$$

for some orthogonal matrix $\mathbf{O}$ and some scalar $d$ that depends on $\mathbf{Z}^{(n)}$. See Randles (2000) for a proof.

The corresponding distances from the origin $d_t^{(n)}(\boldsymbol{\theta}_0, \widehat{\boldsymbol{\Sigma}}^{(n)})$ will be called *pseudo*-Mahalanobis distances, in order to stress the fact that Tyler's estimator of scatter is used instead of the usual sample covariance matrix. The normed standardized residuals $\mathbf{W}_t^{(n)}(\boldsymbol{\theta}_0) := \mathbf{U}_t^{(n)}(\boldsymbol{\theta}_0, \widehat{\boldsymbol{\Sigma}}^{(n)})$—call them *Tyler residuals*—will be used as a multivariate concept of signs.

4.2. *The pseudo-Mahalanobis ranks.* As usual in rank-based nonparametric inference, the pseudo-Mahalanobis distances $d_t^{(n)}(\boldsymbol{\theta}_0, \widehat{\boldsymbol{\Sigma}}^{(n)})$ will be replaced by their ranks. This idea, in the multivariate context, actually goes back to Peters and Randles (1990), who (in a one-sample location context) proved a consistency result, which in the present situation can be stated as follows. Denote by $\hat{R}_t^{(n)}(\boldsymbol{\theta}_0)$ the rank of $d_t^{(n)}(\boldsymbol{\theta}_0, \widehat{\boldsymbol{\Sigma}}^{(n)})$ among $d_1^{(n)}(\boldsymbol{\theta}_0, \widehat{\boldsymbol{\Sigma}}^{(n)}), \ldots, d_n^{(n)}(\boldsymbol{\theta}_0, \widehat{\boldsymbol{\Sigma}}^{(n)})$, and by $R_t^{(n)}(\boldsymbol{\theta}_0, \boldsymbol{\Sigma})$ the rank of $d_t^{(n)}(\boldsymbol{\theta}_0, \boldsymbol{\Sigma})$ among $d_1^{(n)}(\boldsymbol{\theta}_0, \boldsymbol{\Sigma}), \ldots, d_n^{(n)}(\boldsymbol{\theta}_0, \boldsymbol{\Sigma})$.

LEMMA 1 [Peters and Randles (1990)]. *For all $t$, $\hat{R}_t^{(n)}(\boldsymbol{\theta}_0) - R_t^{(n)}(\boldsymbol{\theta}_0, \boldsymbol{\Sigma})$ is $o_{\mathrm{P}}(n)$ as $n \to \infty$ under $\bigcup_f \mathcal{H}^{(n)}(\boldsymbol{\theta}_0, \boldsymbol{\Sigma}, f)$.*

For each $\boldsymbol{\Sigma}$ and $n$, consider the group of *continuous monotone radial transformations* $\mathcal{G}_{\boldsymbol{\Sigma}}^{(n)} = \{\mathcal{G}_g^{(n)}\}$, acting on $(\mathbb{R}^k)^n$, characterized by

$$\begin{aligned}
&\mathcal{G}_g^{(n)}(\mathbf{Z}_1^{(n)}(\boldsymbol{\theta}_0), \ldots, \mathbf{Z}_n^{(n)}(\boldsymbol{\theta}_0)) \\
&\qquad := (g(d_1^{(n)}(\boldsymbol{\theta}_0, \boldsymbol{\Sigma}))\boldsymbol{\Sigma}^{1/2}\mathbf{U}_1^{(n)}(\boldsymbol{\theta}_0, \boldsymbol{\Sigma}), \ldots, g(d_n^{(n)}(\boldsymbol{\theta}_0, \boldsymbol{\Sigma}))\boldsymbol{\Sigma}^{1/2}\mathbf{U}_n^{(n)}(\boldsymbol{\theta}_0, \boldsymbol{\Sigma})),
\end{aligned}$$

where $g : \mathbb{R}^+ \to \mathbb{R}^+$ is continuous, monotone increasing, and such that $g(0) = 0$ and $\lim_{r \to \infty} g(r) = \infty$. The group $\mathcal{G}_{\boldsymbol{\Sigma}}^{(n)}$ is a generating group for the submodel $\bigcup_f \mathcal{H}^{(n)}(\boldsymbol{\theta}_0, \boldsymbol{\Sigma}, f)$, where the union is taken with respect to the set of all possible nonvanishing radial densities. The ranks $R_t^{(n)}(\boldsymbol{\theta}_0, \boldsymbol{\Sigma})$, $t = 1, \ldots, n$, are a maximal invariant for $\mathcal{G}_{\boldsymbol{\Sigma}}^{(n)}$. Lemma 1 thus is an indication that statistics based on the pseudo-Mahalanobis ranks $\hat{R}_t^{(n)}(\boldsymbol{\theta}_0)$ may be asymptotically invariant, in the sense of being asymptotically equivalent to their counterparts based on the unobservable, strictly invariant ranks $R_t^{(n)}(\boldsymbol{\theta}_0, \boldsymbol{\Sigma})$. This will indeed be the case with the test statistics we are proposing (see Propositions 2 and 4).

Note also that the equivariance property (12) of $\mathbf{C}^{(n)}_{\mathrm{Tyl}}$ under affine transformations is sufficient to make the pseudo-Mahalanobis ranks $\hat{R}^{(n)}_t(\boldsymbol{\theta}_0)$ strictly affine-invariant.

4.3. *Tyler residuals.* The transformation $\mathbf{C}^{(n)}_{\mathrm{Tyl}}$ characterized in (11) actually *sphericizes* the problem, in the sense that it transforms elliptically distributed residuals into spherically distributed ones, estimating $\mathbf{U}^{(n)}_t(\boldsymbol{\theta}_0, \boldsymbol{\Sigma})$ by means of the Tyler residuals $\mathbf{W}^{(n)}_t := \mathbf{W}^{(n)}_R(\boldsymbol{\theta}_0) := \mathbf{C}^{(n)}_{\mathrm{Tyl}}\mathbf{Z}^{(n)}_t(\boldsymbol{\theta}_0)/\|\mathbf{C}^{(n)}_{\mathrm{Tyl}}\mathbf{Z}^{(n)}_t(\boldsymbol{\theta}_0)\|$, with the following consistency property.

LEMMA 2. *Under* $\bigcup_f \mathcal{H}^{(n)}(\boldsymbol{\theta}_0, \boldsymbol{\Sigma}, f)$, $\max_{1\le t\le n}\{\|\mathbf{W}^{(n)}_t(\boldsymbol{\theta}_0) - \mathbf{U}^{(n)}_t(\boldsymbol{\theta}_0, \boldsymbol{\Sigma})\|\} = O_{\mathrm{P}}(n^{-1/2})$ *as* $n \to \infty$.

For the proof see the Appendix.

It is clear from (11) that $\mathbf{C}^{(n)}_{\mathrm{Tyl}}(a_1\mathbf{Z}^{(n)}_1, \ldots, a_n\mathbf{Z}^{(n)}_n) = \mathbf{C}^{(n)}_{\mathrm{Tyl}}(\mathbf{Z}^{(n)}_1, \ldots, \mathbf{Z}^{(n)}_n)$ for any real numbers $a_1, \ldots, a_n$, so that $\mathbf{C}^{(n)}_{\mathrm{Tyl}}$ and, therefore, the Tyler residuals $\mathbf{W}^{(n)}_t$ themselves, are strictly invariant under radial monotone transformations. Incidently, it readily follows from (12) that the Tyler residuals enjoy the following strict equivariance property:

LEMMA 3. *Denote by* $\mathbf{W}^{(n)}_t(\mathbf{M})$ *the Tyler residuals computed from the transformed residuals* $\mathbf{M}(\mathbf{Z}^{(n)}_1, \ldots, \mathbf{Z}^{(n)}_n)$. *Then* $\mathbf{W}^{(n)}_t(\mathbf{M}) = \mathbf{O}\mathbf{W}^{(n)}_t$, *where* $\mathbf{O}$ *is the orthogonal matrix in* (12).

For the proof see the Appendix.

Note that Lemma 3 implies that any orthogonally invariant function of the Tyler residuals is strictly affine-invariant. In particular, statistics that are measurable with respect to the cosines of the Euclidean angles between the $\mathbf{W}^{(n)}_t$'s—that is, measurable with respect to the scalar products $(\mathbf{W}^{(n)\prime}_t\mathbf{W}^{(n)}_{\tilde{t}})$—turn out to be affine-invariant. This shows that the Tyler residuals could be used with the same success (consistency, invariance properties) as Randles' interdirections in the construction of the locally asymptotically optimal affine-invariant tests for randomness proposed in Hallin and Paindaveine (2002b). This "angle-based" approach (as opposed to the "interdirection"-based one adopted there) is discussed, for the one-sample location problem, in Hallin and Paindaveine (2002c).

For $k = 1$, the Tyler residuals and pseudo-Mahalanobis ranks reduce to the signs and the ranks of absolute values of the residuals, respectively.

The statistics we are considering in Sections 5 and 6 thus are multivariate generalizations of the serial signed-rank statistics considered in Hallin and Puri (1991).

**5. Rank-based cross-covariance matrices.** The rank-based versions of the cross-covariance matrices (6) we are proposing are of the form

$$\begin{aligned}\widetilde{\mathbf{\Gamma}}^{(n)}_{i;K}(\boldsymbol{\theta}_0) := \mathbf{C}^{(n)\prime}_{\mathrm{Tyl}}\Bigg(&\frac{1}{n-i}\sum_{t=i+1}^{n} K_1\bigg(\frac{\hat{R}^{(n)}_t(\boldsymbol{\theta}_0)}{n+1}\bigg)K_2\bigg(\frac{\hat{R}^{(n)}_{t-i}(\boldsymbol{\theta}_0)}{n+1}\bigg)\\ &\times \mathbf{W}^{(n)}_t(\boldsymbol{\theta}_0)\mathbf{W}^{(n)\prime}_{t-i}(\boldsymbol{\theta}_0)\Bigg)(\mathbf{C}^{(n)\prime}_{\mathrm{Tyl}})^{-1},\end{aligned} \tag{13}$$

where $K_1, K_2 : ]0,1[ \to \mathbb{R}$ are two score functions as in Assumption (C); call (13) a $K$-cross-covariance matrix. Let us shortly review some examples of score functions extending those which are classically considered in univariate rank-based inference. The simplest scores are the constant ones ($K_1(u) = K_2(u) = 1$), and yield multivariate sign cross-covariance matrices

$$\mathbf{C}^{(n)\prime}_{\mathrm{Tyl}}\Bigg(\frac{1}{n-i}\sum_{t=i+1}^{n}\mathbf{W}^{(n)}_t(\boldsymbol{\theta}_0)\mathbf{W}^{(n)\prime}_{t-i}(\boldsymbol{\theta}_0)\Bigg)(\mathbf{C}^{(n)\prime}_{\mathrm{Tyl}})^{-1},$$

leading to serial versions of Randles' multivariate sign test statistic [Randles (2000)]. Linear scores ($K_1(u) = K_2(u) = u$) yield cross-covariance matrices of the Spearman (or Wilcoxon, as only the ranks themselves are involved) type,

$$\begin{aligned}\mathbf{C}^{(n)\prime}_{\mathrm{Tyl}}\Bigg(&\frac{1}{(n-i)\,(n+1)^2}\sum_{t=i+1}^{n}\hat{R}^{(n)}_t(\boldsymbol{\theta}_0)\hat{R}^{(n)}_{t-i}(\boldsymbol{\theta}_0)\\ &\times \mathbf{W}^{(n)}_t(\boldsymbol{\theta}_0)\mathbf{W}^{(n)\prime}_{t-i}(\boldsymbol{\theta}_0)\Bigg)(\mathbf{C}^{(n)\prime}_{\mathrm{Tyl}})^{-1}.\end{aligned} \tag{14}$$

The score functions allowing for local asymptotic optimality under radial density $f_\star$ are $K_1 := \varphi_{f_\star} \circ \tilde{F}^{-1}_{\star k}$ and $K_2 = \tilde{F}^{-1}_{\star k}$ (see Proposition 4). The most familiar example is that of the van der Waerden scores, associated with normal radial densities ($f_\star(r) := \phi(r) = \exp(-r^2/2)$), yielding the van der Waerden cross-covariance matrices

$$\begin{aligned}\mathbf{C}^{(n)\prime}_{\mathrm{Tyl}}\Bigg(&\frac{1}{n-i}\sum_{t=i+1}^{n}\sqrt{\Psi^{-1}_k\bigg(\frac{\hat{R}^{(n)}_t(\boldsymbol{\theta}_0)}{n+1}\bigg)}\\ &\times\sqrt{\Psi^{-1}_k\bigg(\frac{\hat{R}^{(n)}_{t-i}(\boldsymbol{\theta}_0)}{n+1}\bigg)}\mathbf{W}^{(n)}_t(\boldsymbol{\theta}_0)\mathbf{W}^{(n)\prime}_{t-i}(\boldsymbol{\theta}_0)\Bigg)(\mathbf{C}^{(n)\prime}_{\mathrm{Tyl}})^{-1},\end{aligned} \tag{15}$$

where $\Psi_k$ stands for the chi-square distribution function with $k$ degrees of freedom. The Laplace scores, associated with double-exponential radial densities ($f_\star(r) := \exp(-r)$), are another classical example.

In order to study the asymptotic behavior of the $K$-cross-covariance matrices (13) associated with general score functions, under the sequence of null hypotheses as well as under sequences of local alternatives, we first establish the following asymptotic representation and joint normality results; see the Appendix for the proofs.

PROPOSITION 2. *Let Assumptions* (B1) *and* (C) *hold. Then writing*

$$\begin{aligned}\mathbf{\Gamma}^{(n)}_{i;K;\mathbf{\Sigma},f}(\boldsymbol{\theta}_0) &\\ &:= \mathbf{\Sigma}'^{-1/2}\Bigg(\frac{1}{n-i}\sum_{t=i+1}^{n} K_1(\tilde{F}_k(d^{(n)}_t(\boldsymbol{\theta}_0,\mathbf{\Sigma}))) \\ &\qquad\qquad \times K_2(\tilde{F}_k(d^{(n)}_{t-i}(\boldsymbol{\theta}_0,\mathbf{\Sigma})))\mathbf{U}^{(n)}_t(\boldsymbol{\theta}_0,\mathbf{\Sigma})\mathbf{U}^{(n)\prime}_{t-i}(\boldsymbol{\theta}_0,\mathbf{\Sigma})\Bigg)\mathbf{\Sigma}'^{1/2},\end{aligned}$$

$\mathrm{vec}(\widetilde{\mathbf{\Gamma}}^{(n)}_{i;K}(\boldsymbol{\theta}_0) - \mathbf{\Gamma}^{(n)}_{i;K;\mathbf{\Sigma},f}(\boldsymbol{\theta}_0))$ *is* $o_{\mathrm{P}}(n^{-1/2})$ *under* $\mathcal{H}^{(n)}(\boldsymbol{\theta}_0,\mathbf{\Sigma},f)$ *as* $n\to\infty$.

For the proof see the Appendix.

For any square-integrable score function $K$ defined over $]0,1[$, let $\mathrm{E}[K^2(U)] := \int_0^1 K^2(u)\,du$, $D_k(K;f) := \int_0^1 K(u)\,\tilde{F}_k^{-1}(u)\,du$, and $C_k(K;f) := \int_0^1 K(u)\varphi_f\circ\tilde{F}_k^{-1}(u)\,du$. Then we have the following:

PROPOSITION 3. *Let Assumptions* (A), (B1′), (B2) *and* (C) *hold. For any integer* $m$, *the vector*

$$\begin{aligned}\mathbf{S}^{(n)}_{m;K;\mathbf{\Sigma},f}(\boldsymbol{\theta}_0) := ((n-1)^{1/2}(\mathrm{vec}\,\mathbf{\Gamma}^{(n)}_{1;K;\mathbf{\Sigma},f}(\boldsymbol{\theta}_0))', \ldots,& \\ (n-m)^{1/2}(\mathrm{vec}\,\mathbf{\Gamma}^{(n)}_{m;K;\mathbf{\Sigma},f}(\boldsymbol{\theta}_0))')'&\end{aligned} \tag{16}$$

*is asymptotically normal under* $\mathcal{H}^{(n)}(\boldsymbol{\theta}_0,\mathbf{\Sigma},f)$ *and under* $\mathcal{H}^{(n)}(\boldsymbol{\theta}_0+n^{-1/2}\boldsymbol{\tau},\mathbf{\Sigma},f)$, *with mean* $\mathbf{0}$ *under* $\mathcal{H}^{(n)}(\boldsymbol{\theta}_0,\mathbf{\Sigma},f)$ *and mean*

$$\frac{1}{k^2}D_k(K_2;f)C_k(K_1;f)[\mathbf{I}_m\otimes(\mathbf{\Sigma}\otimes\mathbf{\Sigma}^{-1})]\mathbf{Q}^{(m+1)}_{\boldsymbol{\theta}_0}\mathbf{P}_{\boldsymbol{\theta}_0}\mathbf{M}_{\boldsymbol{\theta}_0}\boldsymbol{\tau}$$

*under* $\mathcal{H}^{(n)}(\boldsymbol{\theta}_0+n^{-1/2}\boldsymbol{\tau},\mathbf{\Sigma},f)$, *and with covariance matrix*

$$\frac{1}{k^2}\mathrm{E}[K_1^2(U)]\mathrm{E}[K_2^2(U)][\mathbf{I}_m\otimes(\mathbf{\Sigma}\otimes\mathbf{\Sigma}^{-1})]$$

*under both.*

For the proof see the Appendix.

In order to compare Proposition 3 and the corresponding univariate results in Hallin and Puri (1991, 1994), note that

$$\mathbf{Q}_{\boldsymbol{\theta}_0}^{(m+1)}\mathbf{P}_{\boldsymbol{\theta}_0}\mathbf{M}_{\boldsymbol{\theta}_0}\boldsymbol{\tau}^{(n)} = ((\mathbf{a}_1^{(n)} + \mathbf{b}_1^{(n)})', \ldots, (\mathbf{a}_m^{(n)} + \mathbf{b}_m^{(n)})')',$$

with

$$\mathbf{a}_i^{(n)} := \sum_{j=1}^{\min(p_1,i)} \sum_{k=0}^{i-j} \sum_{l=0}^{\min(q_0,i-j-k)} (\mathbf{G}_{i-j-k-l}\mathbf{B}_l \otimes \mathbf{H}'_k)' \operatorname{vec} \boldsymbol{\gamma}_j^{(n)}$$

and

$$\mathbf{b}_i^{(n)} := \sum_{j=1}^{\min(q_1,i)} (\mathbf{I}_k \otimes \mathbf{H}_{i-j}) \operatorname{vec} \boldsymbol{\delta}_j^{(n)}.$$

Propositions 2 and 3 show that $K$-cross-covariance matrices, while based on multivariate generalizations of signs and ranks, enjoy the same intuitive interpretation and inferential properties as their (traditional) parametric Gaussian counterparts $\boldsymbol{\Gamma}_{i;\boldsymbol{\Sigma},\phi}^{(n)}(\boldsymbol{\theta}_0)$. Proposition 3, for instance, immediately allows for constructing non-Gaussian portmanteau test statistics and deriving their local powers. Just as their classical versions (based on the classical $\boldsymbol{\Gamma}_{i;\boldsymbol{\Sigma},\phi}^{(n)}$'s), such portmanteau tests, however, fail to exploit the information available on the serial dependence structure of the observations, hence, are not optimal. Section 6 is devoted to the construction of locally asymptotically optimal tests based on $K$-cross-covariances.

**6. Optimal tests.** We are now ready to state the main results of this paper: the optimal testing procedures for the problem under study, their invariance and distribution-freeness features, as well as their local powers and optimality properties. Optimality here means local asymptotic minimaxity, either based on fixed-score test statistics, at some selected radial density $f_\star$, or, based on estimated scores, uniformly over some class $\mathcal{F}$ of densities.

6.1. *Fixed-score test statistics.* Letting

$$\begin{aligned}\widetilde{\mathbf{S}}_K^{(n)}(\boldsymbol{\theta}_0) := ((n-1)^{1/2}(\operatorname{vec}\widetilde{\boldsymbol{\Gamma}}_{1;K}^{(n)}(\boldsymbol{\theta}_0))', \ldots,\\ (n-i)^{1/2}(\operatorname{vec}\widetilde{\boldsymbol{\Gamma}}_{i;K}^{(n)}(\boldsymbol{\theta}_0))', \ldots, (\operatorname{vec}\widetilde{\boldsymbol{\Gamma}}_{n-1;K}^{(n)}(\boldsymbol{\theta}_0))')',\end{aligned}$$

define

$$\begin{aligned}n^{1/2}\widetilde{\mathbf{T}}_K^{(n)}(\boldsymbol{\theta}_0) := \mathbf{Q}_{\boldsymbol{\theta}_0}^{(n)\prime}\widetilde{\mathbf{S}}_K^{(n)}(\boldsymbol{\theta}_0) \quad \text{and}\\ \mathbf{J}_{\boldsymbol{\theta}_0,\widehat{\boldsymbol{\Sigma}}}^{(n)} := \mathbf{Q}_{\boldsymbol{\theta}_0}^{(n)\prime}[\mathbf{I}_{n-1} \otimes (\widehat{\boldsymbol{\Sigma}}^{(n)} \otimes \widehat{\boldsymbol{\Sigma}}^{(n)-1})]\mathbf{Q}_{\boldsymbol{\theta}_0}^{(n)},\end{aligned} \tag{17}$$

where $\widehat{\mathbf{\Sigma}}^{(n)}$ denotes Tyler's estimator of scatter (see Section 4.1). Finally, let

$$Q_K^{(n)}(\boldsymbol{\theta}_0) := \frac{k^2 n}{\mathrm{E}[K_1^2(U)]\mathrm{E}[K_2^2(U)]} \widetilde{\mathbf{T}}_K^{(n)\prime}(\boldsymbol{\theta}_0)(\mathbf{J}^{(n)}_{\boldsymbol{\theta}_0,\widehat{\mathbf{\Sigma}}})^{-1}\widetilde{\mathbf{T}}_K^{(n)}(\boldsymbol{\theta}_0).$$

The test statistics $Q_{f_\star}^{(n)}(\boldsymbol{\theta}_0)$ allowing for local asymptotic optimality under radial density $f_\star$ are obtained with the score functions $K_1 := \varphi_{f_\star} \circ \tilde{F}_{\star k}^{-1}$ and $K_2 := \tilde{F}_{\star k}^{-1}$. We then have the following proposition.

Proposition 4. *Assume that Assumptions* (A), (B1), (B2) *and* (C) *hold. Consider the sequence of rank tests $\phi_K^{(n)}$ (resp. $\phi_{f_\star}^{(n)}$) that reject the null hypothesis $\mathcal{H}^{(n)}(\boldsymbol{\theta}_0)$ whenever $Q_K^{(n)}(\boldsymbol{\theta}_0)$ [resp. $Q_{f_\star}^{(n)}(\boldsymbol{\theta}_0)$] exceeds the $\alpha$-upper quantile $\chi^2_{k^2\pi_0,1-\alpha}$ of a chi-square distribution with $k^2\pi_0$ degrees of freedom, where $\pi_0$ is defined in Section 3. Then:*

(i) *the test statistics $Q_K^{(n)}(\boldsymbol{\theta}_0)$ do not depend on the particular choice of the fundamental system $\{\mathbf{\Psi}_t^{(1)}, \ldots, \mathbf{\Psi}_t^{(p_0+q_0)}\}$ (see Section 3); for given values of $p_0$ and $q_0$, they depend on $p_1$ and $q_1$ only through $\pi = \max(p_1 - p_0, q_1 - q_0)$;*

(ii) *$Q_K^{(n)}(\boldsymbol{\theta}_0)$ is asymptotically invariant with respect to the group of continuous monotone radial transformations;*

(iii) *$Q_K^{(n)}(\boldsymbol{\theta}_0)$ is asymptotically chi-square with $k^2\pi_0$ degrees of freedom under $\mathcal{H}^{(n)}(\boldsymbol{\theta}_0)$ (so that $\phi_K^{(n)}$ has asymptotic level $\alpha$), and*

(iv) *$Q_K^{(n)}(\boldsymbol{\theta}_0)$ is asymptotically noncentral chi-square, still with $k^2\pi_0$ degrees of freedom, but with noncentrality parameter*

$$\frac{1}{k^2}\frac{D_k^2(K_2;f)}{\mathrm{E}[K_1^2(U)]}\frac{C_k^2(K_1;f)}{\mathrm{E}[K_2^2(U)]}\boldsymbol{\tau}'\mathbf{N}_{\boldsymbol{\theta}_0,\mathbf{\Sigma}}\boldsymbol{\tau}$$

*under $\mathcal{H}^{(n)}(\boldsymbol{\theta}_0 + n^{-1/2}\boldsymbol{\tau}, \mathbf{\Sigma}, f)$, provided, however, that* (B1) *is reinforced into* (B1′)*, where $\mathbf{N}_{\boldsymbol{\theta}_0,\mathbf{\Sigma}}$ is defined in* (9)*;*

(v) *if we assume that $f_\star$ satifies Assumptions* (B1′), (B2) *and* (C′)*, the sequence of tests $\phi_{f_\star}^{(n)}$ is locally asymptotically maximin for $\mathcal{H}^{(n)}(\boldsymbol{\theta}_0)$ against $\bigcup_{\boldsymbol{\theta}\neq\boldsymbol{\theta}_0}\bigcup_{\mathbf{\Sigma}}\mathcal{H}^{(n)}(\boldsymbol{\theta}, \mathbf{\Sigma}, f_\star)$, at probability level $\alpha$.*

For the proof see the Appendix.

Again, there is no reason to expect the test statistic to be affine-invariant, since the testing problem itself, in general, is not; see Hallin and Paindaveine (2003). Nevertheless, the following proposition establishes that whenever the testing problem under study is affine-invariant (e.g., the problem of testing randomness against VARMA dependence), then the test statistics $Q_K^{(n)}(\boldsymbol{\theta}_0)$ also are affine-invariant.

PROPOSITION 5. (i) *The null hypothesis* $\mathcal{H}^{(n)}(\boldsymbol{\theta}_0)$ *is invariant under affine transformations if and only if* $\boldsymbol{\theta}_0$ *is such that* $\mathbf{A}_i = a_i \mathbf{I}_k$ *for all* $i = 1, \ldots, p_0$ *and* $\mathbf{B}_j = b_j \mathbf{I}_k$ *for all* $j = 1, \ldots, q_0$.

(ii) *When the null hypothesis* $\mathcal{H}^{(n)}(\boldsymbol{\theta}_0)$ *is affine-invariant, then* $Q_K^{(n)}(\boldsymbol{\theta}_0)$ *also is.*

For the proof see the Appendix.

6.2. *Estimated-score test statistics.* The tests $\phi_{f_\star}^{(n)}$ considered in Proposition 4 achieve parametric efficiency at radial density $f_\star$. ARMA models, though, under adequate assumptions, are adaptive; this has been shown formally in the univariate case only [without even requiring symmetric innovation densities; see, e.g., Drost, Klaassen and Werker (1997)], but is very likely to hold also in higher dimensions. *Adaptive* optimality property—that is, parametric optimality at all $f_\star$—thus can be expected, provided that *estimated scores* are considered. Proposition 6 shows that this, indeed, is the case.

An adaptive procedure could be based on the score function $\varphi_{\hat{f}}$ associated with an adequate estimator $\hat{f}$ of the radial density. While being uniformly locally asymptotically maximin, such a procedure, however, would not have the very desirable properties of rank-based procedures. This is why we rather propose, in the spirit of Hallin and Werker (2003), an adaptive version of the rank-based procedures described in Proposition 4.

Let us first assume that $\boldsymbol{\Sigma}$ is known, so that the genuine distances $d_t^{(n)} := d_t^{(n)}(\boldsymbol{\theta}_0, \boldsymbol{\Sigma})$ can be computed from the observations. Denote by $R_t^{(n)} := R_t^{(n)}(\boldsymbol{\theta}_0, \boldsymbol{\Sigma})$ the rank of $d_t^{(n)}$ among $d_1^{(n)}, \ldots, d_n^{(n)}$: under $\mathcal{H}^{(n)}(\boldsymbol{\theta}_0, \boldsymbol{\Sigma}, f)$ the $R_t^{(n)}$'s are the ranks of i.i.d. random variables with probability density function $\tilde{f}_k$. Next consider any continuous kernel density estimator $\check{f}_k^{(n)}$ of $\tilde{f}_k$ that is measurable with respect to the order statistic of the $d_t^{(n)}$'s and satisfies

$$\begin{aligned}
\mathrm{E}\Bigg[\Bigg[&\varphi_{\check{f}_k^{(n)}}\left((\check{F}_k^{(n)})^{-1}\left(\frac{R_t^{(n)}}{n+1}\right)\right)(\check{F}_k^{(n)})^{-1}\left(\frac{R_{t-i}^{(n)}}{n+1}\right) \\
&- \varphi_{\tilde{f}_k}\left(\tilde{F}_k^{-1}\left(\frac{R_t^{(n)}}{n+1}\right)\right)\tilde{F}_k^{-1}\left(\frac{R_{t-i}^{(n)}}{n+1}\right)\Bigg]^2 \Bigg| \check{f}_k^{(n)}\Bigg] = o_{\mathrm{P}}(1),
\end{aligned} \tag{18}$$

under $\mathcal{H}^{(n)}(\boldsymbol{\theta}_0, \boldsymbol{\Sigma}, f)$ as $n \to \infty$, where $\check{F}_k^{(n)}$ denotes the cumulative distribution function associated with $\check{f}_k^{(n)}$.

A possible choice for $\check{f}_k^{(n)}$ satisfying (18) is given in Hájek and Šidák [(1967), (1.5.7) of Chapter VII]. Another one, specifically constructed for radial densities, is proposed in Liebscher (2005). An adaptive (still, under specified $\mathbf{\Sigma}$) version of (13) is then

$$\check{\mathbf{\Gamma}}^{(n)}_{i;\mathbf{\Sigma}}(\boldsymbol{\theta}_0) := \mathbf{\Sigma}'^{-1/2}\Bigg(\frac{1}{n-i}\sum_{t=i+1}^{n}\check{\varphi}^{(n)}_f\Bigg((\check{F}^{(n)}_k)^{-1}\bigg(\frac{R^{(n)}_t}{n+1}\bigg)\Bigg) \times (\check{F}^{(n)}_k)^{-1}\bigg(\frac{R^{(n)}_{t-i}}{n+1}\bigg) \times \mathbf{U}^{(n)}_t(\boldsymbol{\theta}_0,\mathbf{\Sigma})\mathbf{U}^{(n)\prime}_{t-i}(\boldsymbol{\theta}_0,\mathbf{\Sigma})\Bigg)\mathbf{\Sigma}'^{1/2}, \tag{19}$$

where we let $\check{\varphi}^{(n)}_f(r) := \varphi_{\check{f}^{(n)}_k}(r) + (k-1)/r$ [since $\varphi_f(r) = \varphi_{\tilde{f}_k}(r) + (k-1)/r$].

Of course, in practice $\mathbf{\Sigma}$ is not known, and only the estimated distances $\hat{d}^{(n)}_t := d^{(n)}_t(\boldsymbol{\theta}_0, \hat{\mathbf{\Sigma}}^{(n)})$ can be computed: instead of $\check{\mathbf{\Gamma}}^{(n)}_{i;\mathbf{\Sigma}}(\boldsymbol{\theta}_0)$ given in (19), we therefore rather use (with the notation of Section 4)

$$\hat{\mathbf{\Gamma}}^{(n)}_{i}(\boldsymbol{\theta}_0) := \mathbf{C}^{(n)\prime}_{\mathrm{Tyl}}\Bigg(\frac{1}{n-i}\sum_{t=i+1}^{n}\hat{\varphi}^{(n)}_f\Bigg((\hat{F}^{(n)}_k)^{-1}\bigg(\frac{\hat{R}^{(n)}_t}{n+1}\bigg)\Bigg) \times (\hat{F}^{(n)}_k)^{-1}\bigg(\frac{\hat{R}^{(n)}_{t-i}}{n+1}\bigg) \times \mathbf{W}^{(n)}_t(\boldsymbol{\theta}_0)\mathbf{W}^{(n)\prime}_{t-i}(\boldsymbol{\theta}_0)\Bigg)(\mathbf{C}^{(n)\prime}_{\mathrm{Tyl}})^{-1}, \tag{20}$$

where $\check{f}^{(n)}_k$, $\check{F}^{(n)}_k$ and $\check{\varphi}^{(n)}_f$ have been replaced by their counterparts $\hat{f}^{(n)}_k$, $\hat{F}^{(n)}_k$ and $\hat{\varphi}^{(n)}_f$ computed from the order statistic of the $\hat{d}^{(n)}_t$'s. Using the multivariate Slutsky theorem and working as in the proof of Proposition 2, we obtain that the difference between (19) and (20) is $o_{\mathrm{P}}(n^{-1/2})$ under $\bigcup_f \mathcal{H}^{(n)}(\boldsymbol{\theta}_0, \mathbf{\Sigma}, f)$ as $n \to \infty$. A direct adaptation of the proof of Proposition 3.4 in Hallin and Werker (2003) then yields a multivariate generalization of the (symmetric version of) Proposition 6.4 in Hallin and Werker (1999). This adaptation, however, requires the Fisher information for location associated with $\tilde{f}_k$ to be finite. Denote by $\mathcal{F}$ the set of all radial densities $f$ for which this condition is satisfied: clearly, $\{f | \mathcal{I}_{k,f} < \infty$ and $\int_0^\infty r^{k-3} f(r)\,dr < \infty\} \subset \mathcal{F}$ and, in the univariate case ($k = 1$), $\mathcal{F} = \{f | \mathcal{I}_{1,f} < \infty\}$.

LEMMA 4. *Let Assumptions* (B1) *and* (B2) *hold, and assume that* $f \in \mathcal{F}$ *satisfies Assumption* (C′). *Then both* $\mathrm{vec}(\check{\mathbf{\Gamma}}^{(n)}_i(\boldsymbol{\theta}_0) - \mathbf{\Gamma}^{(n)}_{i;\mathbf{\Sigma},f}(\boldsymbol{\theta}_0))$ *and* $\mathrm{vec}(\hat{\mathbf{\Gamma}}^{(n)}_i(\boldsymbol{\theta}_0) - \mathbf{\Gamma}^{(n)}_{i;\mathbf{\Sigma},f}(\boldsymbol{\theta}_0))$ *are* $o_{\mathrm{P}}(n^{-1/2})$ *under* $\mathcal{H}^{(n)}(\boldsymbol{\theta}_0, \mathbf{\Sigma}, f)$ *as* $n \to \infty$.

In order to construct adaptive procedures, we still need to estimate the asymptotic variance–covariance matrices of either (19) or (20). More precisely, we need consistent estimates of $\mathcal{I}_{k,f}$ and $v_{k,f} := \mu_{k+1;f}/\mu_{k-1;f} =$

$\mathrm{E}[(\tilde{F}_k^{-1}(U))^2]$. Such estimates are provided by

$$\check{\mathcal{I}}^{(n)} := \frac{1}{n}\sum_{t=1}^{n}\biggl(\check{\varphi}_f^{(n)}\circ(\check{F}_k^{(n)})^{-1}\biggl(\frac{R_t^{(n)}}{n+1}\biggr)\biggr)^2$$

and

$$\check{v}^{(n)} := \frac{1}{n}\sum_{t=1}^{n}\biggl((\check{F}_k^{(n)})^{-1}\biggl(\frac{R_t^{(n)}}{n+1}\biggr)\biggr)^2,$$

$$\hat{\mathcal{I}}^{(n)} := \frac{1}{n}\sum_{t=1}^{n}\biggl(\hat{\varphi}_f^{(n)}\circ(\hat{F}_k^{(n)})^{-1}\biggl(\frac{\hat{R}_t^{(n)}}{n+1}\biggr)\biggr)^2$$

and

$$\hat{v}^{(n)} := \frac{1}{n}\sum_{t=1}^{n}\biggl((\hat{F}_k^{(n)})^{-1}\biggl(\frac{\hat{R}_t^{(n)}}{n+1}\biggr)\biggr)^2,$$

respectively. Note that the product $\check{\mathcal{I}}^{(n)}\check{v}^{(n)}$ (resp. $\hat{\mathcal{I}}^{(n)}\hat{v}^{(n)}$) depends on the estimated radial density $\check{f}_k^{(n)}$ (resp. $\hat{f}_k^{(n)}$) only through its *density type*—namely, the scale family $\{a\check{f}_k^{(n)}(ar), a>0\}$ [resp. $\{a\hat{f}_k^{(n)}(ar), a>0\}$]. By the way, the same property holds true for the adaptive rank-based cross-covariances (19) and (20); consequently, without any loss of generality, we may assume that $\check{f}_k^{(n)}$ and $\hat{f}_k^{(n)}$ are such that $\check{v}^{(n)} = \hat{v}^{(n)} = 1$.

Defining

$$\check{\mathbf{S}}^{(n)}(\boldsymbol{\theta}_0) := \check{\mathbf{S}}_{\boldsymbol{\Sigma}}^{(n)}(\boldsymbol{\theta}_0) := ((n-1)^{1/2}(\operatorname{vec}\check{\boldsymbol{\Gamma}}_{1;\boldsymbol{\Sigma}}^{(n)}(\boldsymbol{\theta}_0))', \ldots, (\operatorname{vec}\check{\boldsymbol{\Gamma}}_{n-1;\boldsymbol{\Sigma}}^{(n)}(\boldsymbol{\theta}_0))')'$$

and

$$n^{1/2}\check{\mathbf{T}}^{(n)}(\boldsymbol{\theta}_0) := n^{1/2}\check{\mathbf{T}}_{\boldsymbol{\Sigma}}^{(n)}(\boldsymbol{\theta}_0) := \mathbf{Q}_{\boldsymbol{\theta}_0}^{(n)\prime}\check{\mathbf{S}}^{(n)}(\boldsymbol{\theta}_0),$$

let

$$\check{Q}^{(n)}(\boldsymbol{\theta}_0) := \check{Q}_{\boldsymbol{\Sigma}}^{(n)}(\boldsymbol{\theta}_0) := \frac{k^2 n}{\check{\mathcal{I}}^{(n)}\check{v}^{(n)}}\check{\mathbf{T}}^{(n)\prime}(\boldsymbol{\theta}_0)(\mathbf{J}_{\boldsymbol{\theta}_0,\boldsymbol{\Sigma}}^{(n)})^{-1}\check{\mathbf{T}}^{(n)}(\boldsymbol{\theta}_0), \tag{21}$$

where $\mathbf{J}_{\boldsymbol{\theta}_0,\boldsymbol{\Sigma}}^{(n)}$ is defined in (17). The same quantities, when computed from the $\hat{\boldsymbol{\Gamma}}_i^{(n)}(\boldsymbol{\theta}_0)$'s, are denoted by $\hat{\mathbf{S}}^{(n)}(\boldsymbol{\theta}_0)$ and $\hat{\mathbf{T}}^{(n)}(\boldsymbol{\theta}_0)$, respectively, yielding the test statistic

$$\hat{Q}^{(n)}(\boldsymbol{\theta}_0) := \frac{k^2 n}{\hat{\mathcal{I}}^{(n)}\hat{v}^{(n)}}\,\hat{\mathbf{T}}^{(n)\prime}(\boldsymbol{\theta}_0)(\mathbf{J}_{\boldsymbol{\theta}_0,\hat{\boldsymbol{\Sigma}}}^{(n)})^{-1}\hat{\mathbf{T}}^{(n)}(\boldsymbol{\theta}_0) = \hat{Q}_{\hat{\boldsymbol{\Sigma}}}^{(n)}(\boldsymbol{\theta}_0).$$

The test statistic (21) has the very desirable property of being conditionally distribution-free. Conditional upon the $\sigma$-algebra $\mathcal{D}^{(n)}$ generated by the order statistic $\mathbf{d}_{(\cdot)}^{(n)}$ of the *exact* distances $\mathbf{d}^{(n)} := (d_1^{(n)}, \ldots, d_n^{(n)})$, indeed:

(a) the vector of ranks $\mathbf{R}^{(n)} := (R_1^{(n)}, \dots, R_n^{(n)})$ is uniformly distributed over the $n!$ permutations of $(1, \dots, n)$,

(b) the normalized residuals $\mathbf{U}_t^{(n)}$ are i.i.d. and uniformly distributed over the unit hypersphere, and

(c) the ranks $\mathbf{R}^{(n)}$ and the residuals $\mathbf{U}_t^{(n)}$ are mutually independent.

The situation is thus entirely parallel to the classical case of univariate signed ranks: conditional on $\mathcal{D}^{(n)}$, $\check{Q}^{(n)}(\boldsymbol{\theta}_0)$ is distribution-free. Denote by $\check{q}_\alpha(\mathbf{d}_{(\cdot)}^{(n)})$ its upper $\alpha$-quantile, and by $\check{\phi}^{(n)} := \check{\phi}_{\mathbf{\Sigma}}^{(n)}$ the indicator of the event $\check{Q}^{(n)}(\boldsymbol{\theta}_0) > \check{q}_\alpha(\mathbf{d}_{(\cdot)}^{(n)})$. This test actually has Neyman $\alpha$-structure with respect to $\mathbf{d}_{(\cdot)}^{(n)}$ and, consequently, is a permutation test. Proposition 6 and Lemma 4, moreover, imply that the sequence $\check{\phi}^{(n)}$ is asymptotically optimal, uniformly in $f$, against $\bigcup_f \mathcal{H}^{(n)}(\boldsymbol{\theta}_0, \mathbf{\Sigma}, f)$.

Unfortunately, unlike univariate adaptive signed rank tests, this permutation test cannot be implemented, since $\mathbf{\Sigma}$, in practice, is unspecified. Instead of $\check{\phi}^{(n)}$, based on $\check{Q}^{(n)}(\boldsymbol{\theta}_0)$, we therefore recommend $\hat{\phi}^{(n)} = \check{\phi}_{\hat{\mathbf{\Sigma}}}^{(n)}$, based on the test statistic $\hat{Q}^{(n)}(\boldsymbol{\theta}_0)$, which rejects the null hypothesis $\mathcal{H}^{(n)}(\boldsymbol{\theta}_0)$ whenever $\hat{Q}^{(n)}(\boldsymbol{\theta}_0)$ exceeds the $\alpha$-upper quantile $\chi^2_{k^2\pi_0, 1-\alpha}$ of a chi-square distribution with $k^2\pi_0$ degrees of freedom. In view of Lemma 4, $\check{Q}^{(n)}(\boldsymbol{\theta}_0)$ and $\hat{Q}^{(n)}(\boldsymbol{\theta}_0)$ are asymptotically equivalent under $\bigcup_f \mathcal{H}^{(n)}(\boldsymbol{\theta}_0, \mathbf{\Sigma}, f)$ and contiguous alternatives: $\hat{\phi}^{(n)}$ and $\check{\phi}^{(n)}$ thus share the same asymptotic optimality properties. On the other hand, $\hat{\phi}^{(n)}$ loses the attractive finite-sample Neyman $\alpha$-structure of $\check{\phi}^{(n)}$.

Summing up, the following proposition is a direct consequence of Lemma 4.

Proposition 6. *Let Assumptions* (A), (B1′) *and* (B2) *hold, and assume that $f \in \mathcal{F}$ satisfies Assumption* (C′). *Then:*

(i) *statements* (i)–(iii) *of Proposition 4 hold for $\hat{\phi}^{(n)}$; statement* (iv) *also holds, with asymptotic noncentrality parameter*

$$\frac{1}{k^2}\mathrm{E}[(\tilde{F}_k^{-1}(U))^2]\mathrm{E}[(\varphi_f(\tilde{F}_k^{-1}(U)))^2]\boldsymbol{\tau}'\mathbf{N}_{\boldsymbol{\theta}_0,\mathbf{\Sigma}}\boldsymbol{\tau}$$

*under $\mathcal{H}^{(n)}(\boldsymbol{\theta}_0 + n^{-1/2}\boldsymbol{\tau}, \mathbf{\Sigma}, f)$;*

(ii) *the sequence of tests $\hat{\phi}^{(n)}$ is locally asymptotically maximin for $\mathcal{H}^{(n)}(\boldsymbol{\theta}_0)$ against $\bigcup_{\boldsymbol{\theta} \neq \boldsymbol{\theta}_0} \bigcup_{\mathbf{\Sigma}} \bigcup_f \mathcal{H}^{(n)}(\boldsymbol{\theta}, \mathbf{\Sigma}, f)$, at probability level $\alpha$, where the third union is taken over all radial densities $f \in \mathcal{F}$ satisfying Assumptions* (B1′), (B2) *and* (C′).

Proposition 5 readily extends to this adaptive procedure.

6.3. *The Gaussian procedure.* We now briefly describe the parametric Gaussian procedure for the problem treated in Propositions 4 and 6. This Gaussian test will serve as a benchmark in Section 7 for the computation of asymptotic relative efficiencies.

Under Gaussian assumptions the empirical covariance $\mathbf{S}^{(n)} := n^{-1} \times \sum_{t=1}^{n} \mathbf{Z}_t^{(n)}(\boldsymbol{\theta}_0)\mathbf{Z}_t^{(n)\prime}(\boldsymbol{\theta}_0)$ is a consistent estimator under $\mathcal{H}^{(n)}(\boldsymbol{\theta}_0, \boldsymbol{\Sigma}, f)$ of the innovation covariance $(\mathrm{E}[(\tilde{F}_k^{-1}(U))^2]/k)\boldsymbol{\Sigma}$. Let

$$\mathbf{J}_{\mathcal{N};\boldsymbol{\theta}_0}^{(n)} := \mathbf{Q}_{\boldsymbol{\theta}_0}^{(n)\prime}[\mathbf{I}_{n-1} \otimes \hat{\boldsymbol{\Gamma}}_{\boldsymbol{\theta}_0}^{(n)}]\mathbf{Q}_{\boldsymbol{\theta}_0}^{(n)},$$

where $\hat{\boldsymbol{\Gamma}}_{\boldsymbol{\theta}_0}^{(n)} := (n-1)^{-1} \sum_{t=2}^{n} \operatorname{vec}(\mathbf{Z}_t^{(n)}(\boldsymbol{\theta}_0)\mathbf{Z}_{t-1}^{(n)\prime}(\boldsymbol{\theta}_0))(\operatorname{vec}(\mathbf{Z}_t^{(n)}(\boldsymbol{\theta}_0)\mathbf{Z}_{t-1}^{(n)\prime}(\boldsymbol{\theta}_0)))'$. In view of the ergodic theorem [see Hannan (1970), Theorem 2, page 203], $\hat{\boldsymbol{\Gamma}}_{\boldsymbol{\theta}_0}^{(n)}$ is consistent under $\mathcal{H}^{(n)}(\boldsymbol{\theta}_0, \boldsymbol{\Sigma}, f)$ for $(\mathrm{E}[(\tilde{F}_k^{-1}(U))^2]/k)^2 \boldsymbol{\Sigma} \otimes \boldsymbol{\Sigma}^{-1}$. The following proposition then follows along the same lines as Proposition 4.

PROPOSITION 7. *Let Assumptions* (A), (B1′) *and* (B2) *hold. Define*

$$Q_{\mathcal{N}}^{(n)}(\boldsymbol{\theta}_0) := n\, \mathbf{T}_{\mathbf{S},\phi}^{(n)\prime}(\boldsymbol{\theta}_0)(\mathbf{J}_{\mathcal{N};\boldsymbol{\theta}_0}^{(n)})^{-1}\mathbf{T}_{\mathbf{S},\phi}^{(n)}(\boldsymbol{\theta}_0). \tag{22}$$

*Consider the sequence of parametric Gaussian tests $\phi_{\mathcal{N}}^{(n)}$ rejecting the null hypothesis $\mathcal{H}^{(n)}(\boldsymbol{\theta}_0)$ whenever $Q_{\mathcal{N}}^{(n)}(\boldsymbol{\theta}_0)$ exceeds the $\alpha$-upper quantile $\chi^2_{k^2\pi_0,1-\alpha}$ of a chi-square distribution with $k^2\pi_0$ degrees of freedom. Then*:

(i) *statements* (i) *and* (iii) *in Proposition 4 hold for $\phi_{\mathcal{N}}^{(n)}$; statement* (iv) *also holds, with asymptotic noncentrality parameter* $(\mathrm{E}^2[\tilde{F}_k^{-1}(U)\varphi_f(\tilde{F}_k^{-1}(U))]/k^2) \times \boldsymbol{\tau}'\mathbf{N}_{\boldsymbol{\theta}_0,\boldsymbol{\Sigma}}\boldsymbol{\tau}$ *under* $\mathcal{H}^{(n)}(\boldsymbol{\theta}_0 + n^{-1/2}\boldsymbol{\tau}, \boldsymbol{\Sigma}, f)$;

(ii) *the sequence of tests $\phi_{\mathcal{N}}^{(n)}$ is locally asymptotically maximin for $\mathcal{H}^{(n)}(\boldsymbol{\theta}_0)$ against the Gaussian alternative $\bigcup_{\boldsymbol{\Sigma}} \mathcal{H}^{(n)}(\boldsymbol{\theta}_0, \boldsymbol{\Sigma}, \phi)$, at probability level $\alpha$.*

The test statistic $Q_{\mathcal{N}}^{(n)}(\boldsymbol{\theta}_0)$ is not (not even asymptotically) invariant under continuous monotone radial transformations. However, it is asymptotically distribution-free. On the other hand, $Q_{\mathcal{N}}^{(n)}(\boldsymbol{\theta}_0)$, just like $Q_K^{(n)}(\boldsymbol{\theta}_0)$ and $\hat{Q}^{(n)}(\boldsymbol{\theta}_0)$, is affine-invariant whenever the null hypothesis is.

## 7. Asymptotic performance.

7.1. *Asymptotic relative efficiencies.* Computing the ratios of the noncentrality parameters in the asymptotic distributions of $Q_K^{(n)}$, $Q_{f_\star}^{(n)}$ and $\hat{Q}^{(n)}$ with respect to $Q_{\mathcal{N}}^{(n)}$ (see Propositions 4, 6 and 7) yields the asymptotic relative efficiencies of these tests with respect to their parametric Gaussian counterparts.

PROPOSITION 8. *Let Assumptions* (A), (B1′), (B2) *and* (C) *hold. Then*:

(i) *the asymptotic relative efficiency under radial density* $f$ *of* $\phi_K^{(n)}$ *with respect to* $\phi_{\mathcal{N}}^{(n)}$ *is*

$$\mathrm{ARE}_{k,f}(\phi_K^{(n)}/\phi_{\mathcal{N}}^{(n)}) = \frac{1}{k^2}\frac{D_k^2(K_2; f)}{\mathrm{E}[K_2^2(U)]}\frac{C_k^2(K_1; f)}{\mathrm{E}[K_1^2(U)]};$$

(ii) *assuming that* (C′) *holds instead of* (C),

$$\mathrm{ARE}_{k,f}(\phi_{f_\star}^{(n)}/\phi_{\mathcal{N}}^{(n)}) = \frac{1}{k^2}\frac{D_k^2(f_\star, f)}{D_k(f_\star)}\frac{C_k^2(f_\star, f)}{C_k(f_\star)},$$

*where we write* $D_k(f_1, f_2)$ *and* $C_k(f_1, f_2)$ *for* $D_k(\tilde{F}_{1k}^{-1}; f_2)$ *and* $C_k(\varphi_{f_1} \circ \tilde{F}_{1k}^{-1}; f_2)$, *respectively, and let* $C_k(f_\star) := C_k(f_\star, f_\star)$ *and* $D_k(f_\star) := D_k(f_\star, f_\star)$;

(iii) *assuming, moreover, that* $f \in \mathcal{F}$ *satisfies Assumption* (C′), *the asymptotic relative efficiency of the adaptive test* $\hat{\phi}^{(n)}$ *with respect to* $\phi_{\mathcal{N}}^{(n)}$ *under the radial density* $f$ *is*

$$\mathrm{ARE}_{k,f}(\hat{\phi}^{(n)}/\phi_{\mathcal{N}}^{(n)}) = \frac{1}{k^2} D_k(f) C_k(f).$$

The AREs for the fixed-score procedures obtained in Proposition 8 coincide with those obtained in Hallin and Paindaveine (2002b) for the related problem of testing randomness against VARMA dependence. The numerical values of AREs of several versions of the proposed procedures (van der Waerden and Laplace score tests, sign test, Spearman-type test) with respect to the Gaussian procedure, under multivariate $t$-distributions with various degrees of freedom, are reported there. As usual in rank-based inference, the gain of efficiency over parametric $L^2$ procedures increases with the tail weight [see Hallin and Paindaveine (2002b)].

In this section, we thus concentrate on the adaptive procedure described in Proposition 6. As in Randles (1989), consider the family of power-exponential distributions with density

$$\underline{f_\nu}(\mathbf{x}) = K_{k,\nu}\frac{1}{(\det \boldsymbol{\Sigma})^{1/2}} \exp[-((\mathbf{x}-\boldsymbol{\theta})'\boldsymbol{\Sigma}^{-1}(\mathbf{x}-\boldsymbol{\theta})/c_0)^\nu], \qquad \nu > 0, \tag{23}$$

with

$$c_0 := \frac{k\,\Gamma(k/2\nu)}{\Gamma((k+2)/2\nu)} \quad \text{and} \quad K_{k,\nu} := \frac{\nu\Gamma(k/2)}{\Gamma(k/2\nu)(\pi c_0)^{k/2}}.$$

This family corresponds to radial densities of the form $f_\nu(r) := \exp[-(r^2/c_0)^\nu]$, and allows for considering a variety of tail weights indexed by $\nu$. The $k$-variate normal case corresponds to $\nu = 1$, while, for $0 < \nu < 1$ (resp. $\nu > 1$), the tails are heavier (resp. lighter) than in the normal case.

TABLE 1
*Asymptotic relative efficiencies of the adaptive test* $\hat{\phi}^{(n)}$ *w.r.t. the Gaussian test* $\phi_{\mathcal{N}}^{(n)}$ *in the elliptically symmetric power-exponential family* (23), *for various values of the tail index* $\nu$ *and the space dimension* $k$

| | $\boldsymbol{\nu}$ | | | | | | | |
|---|---|---|---|---|---|---|---|---|
| $\boldsymbol{k}$ | **0.1** | **0.2** | **0.3** | **0.5** | **1** | **2** | **5** | **10** |
| 1 | — | — | 28.40 | 2.00 | 1.00 | 1.37 | 3.18 | 6.43 |
| 3 | 261.24 | 8.08 | 2.77 | 1.33 | 1.00 | 1.22 | 2.30 | 4.26 |
| 4 | 59.63 | 4.77 | 2.16 | 1.25 | 1.00 | 1.18 | 2.08 | 3.71 |
| 6 | 14.81 | 2.84 | 1.69 | 1.17 | 1.00 | 1.13 | 1.81 | 3.03 |
| 8 | 7.51 | 2.19 | 1.48 | 1.13 | 1.00 | 1.10 | 1.65 | 2.63 |
| 10 | 5.02 | 1.88 | 1.37 | 1.10 | 1.00 | 1.09 | 1.54 | 2.36 |
| $\infty$ | 1.00 | 1.00 | 1.00 | 1.00 | 1.00 | 1.00 | 1.00 | 1.00 |

Provided that $4\nu + k - 2 > 0$, Proposition 8 yields

$$\text{ARE}_{k,f_\nu}(\hat{\phi}^{(n)}/\phi_{\mathcal{N}}^{(n)}) = \frac{4\nu^2}{k^2}\frac{\Gamma((k+2)/2\nu)\Gamma((4\nu+k-2)/2\nu)}{\Gamma^2(k/2\nu)}. \tag{24}$$

Table 1 above provides some numerical values of (24).

7.2. *A multivariate version of two classical univariate results.* Since the AREs obtained in Proposition 8 for the fixed-score procedures $\phi_K^{(n)}$ and $\phi_{f_\star}^{(n)}$ coincide with those in Hallin and Paindaveine (2002b), the generalizations obtained there of the famous Chernoff–Savage and Hodges–Lehmann results still hold here. In view of their importance, we adapt these results to the present context, referring to Hallin and Paindaveine (2002b) for proofs and details.

*A multivariate serial Chernoff–Savage result.* As in the univariate case, the van der Waerden version of the proposed rank-based procedure is uniformly more efficient than the corresponding parametric Gaussian procedure. More precisely, the following generalization of the results of Chernoff and Savage (1958) and Hallin (1994) holds.

PROPOSITION 9. *Let Assumption* (A) *hold. Denote by* $\phi_{\text{vdW}}^{(n)}$ *and* $\phi_{\mathcal{N}}^{(n)}$ *the van der Waerden test, based on the cross-covariance matrices* (15), *and the Gaussian test based on the test statistic* (22), *respectively. For any* $f$ *satisfying Assumptions* (B1′) *and* (B2),

$$\text{ARE}_{k,f}(\phi_{\text{vdW}}^{(n)}/\phi_{\mathcal{N}}^{(n)}) \geq 1,$$

*where equality holds if and only if* $f$ *is normal.*

TABLE 2
*Some numerical values, for various values $k$ of the space dimension, of the lower bound for the asymptotic relative efficiency of the Spearman test $\phi_{\mathrm{SP}}^{(n)}$ with respect to the Gaussian one $\phi_{\mathcal{N}}^{(n)}$*

| $k$ | $\inf_f \mathrm{ARE}_{k,f}(\phi_{\mathrm{SP}}^{(n)}/\phi_{\mathcal{N}}^{(n)})$ | $k$ | $\inf_f \mathrm{ARE}_{k,f}(\phi_{\mathrm{SP}}^{(n)}/\phi_{\mathcal{N}}^{(n)})$ |
|---|---|---|---|
| 1 | 0.856 | 5 | 0.818 |
| 2 | 0.913 | 6 | 0.797 |
| 3 | 0.878 | 10 | 0.742 |
| 4 | 0.845 | $+\infty$ | 0.563 |

*A multivariate serial Hodges–Lehmann result.* Denote by $Q_{\mathrm{SP}}^{(n)}(\boldsymbol{\theta}_0)$ the Spearman-type version of the test statistics $Q_K^{(n)}(\boldsymbol{\theta}_0)$, based on the cross-covariances (14) associated with linear scores. This statistic can be considered as the angle-based serial version of Peters and Randles' Wilcoxon-type test statistic [see Hallin and Paindaveine (2002c) and Peters and Randles (1990)].

Although the resulting test $\phi_{\mathrm{SP}}^{(n)}$ is never optimal [there is no $f_\star$ such that $Q_{f_\star}^{(n)}(\boldsymbol{\theta}_0)$ coincides with $Q_{\mathrm{SP}}^{(n)}(\boldsymbol{\theta}_0)$], the resulting Spearman-type procedure exhibits excellent asymptotic efficiency properties, especially for relatively small dimensions $k$. To show this, we extend Hodges and Lehmann's (1956) celebrated "0.864 result" by computing, for any dimension $k$, the lower bound for the asymptotic relative efficiency of $\phi_{\mathrm{SP}}^{(n)}$ with respect to the Gaussian procedure $\phi_{\mathcal{N}}^{(n)}$. More precisely, we have the following proposition [see Hallin and Paindaveine (2002b) for the proof].

PROPOSITION 10. *Let Assumption* (A) *hold. Define*

$$c_k := \inf\{x > 0 | (\sqrt{x} J_{\sqrt{2k-1}/2}(x))' = 0\},$$

*where $J_r$ denotes the Bessel function of the first kind of order $r$. The lower bound for the asymptotic relative efficiency of $\phi_{\mathrm{SP}}^{(n)}$ with respect to $\phi_{\mathcal{N}}^{(n)}$ is*

$$\inf_f \mathrm{ARE}_{k,f}(\phi_{\mathrm{SP}}^{(n)}/\phi_{\mathcal{N}}^{(n)}) = 9(2c_k^2 + k - 1)^4/2^{10}k^2c_k^4, \tag{25}$$

*where the infimum is taken over all radial densities $f$ satisfying Assumptions* (B1′) *and* (B2).

Some numerical values are given in Table 2. Note that the sequence of lower bounds (25) is monotonically decreasing in $k$ for $k \geq 2$, and tends to $9/16 = 0.5625$ as $k \to \infty$.

# APPENDIX

## A.1. Proofs of Proposition 1 and Lemmas 2 and 3.

PROOF OF PROPOSITION 1. Garel and Hallin (1995) show that the linear part in the quadratic approximation of $L^{(n)}_{\boldsymbol{\theta}_0+n^{-1/2}\boldsymbol{\tau}^{(n)}/\boldsymbol{\theta}_0;\boldsymbol{\Sigma},f}$ can be written as

$$\boldsymbol{\tau}^{(n)\prime}\boldsymbol{\Delta}^{(n)}_{\boldsymbol{\Sigma},f}(\boldsymbol{\theta}_0)=\sum_{i=1}^{n-1}(n-i)^{1/2}\operatorname{tr}[\mathbf{d}_i^{(n)\prime}(\boldsymbol{\theta}_0)\boldsymbol{\Gamma}^{(n)}_{i;\boldsymbol{\Sigma},f}(\boldsymbol{\theta}_0)],$$

where

$$\mathbf{d}_i^{(n)}(\boldsymbol{\theta}_0):=\sum_{j=1}^{\min(p_1,i)}\sum_{k=0}^{i-j}\sum_{l=0}^{\min(q_0,i-j-k)}\mathbf{H}_k\gamma_j^{(n)}\mathbf{G}_{i-j-k-l}\mathbf{B}_l+\sum_{j=1}^{\min(q_1,i)}\mathbf{H}_{i-j}\boldsymbol{\delta}_j^{(n)}.$$

Using $\operatorname{tr}(\mathbf{AB})=(\operatorname{vec}\mathbf{A}')'(\operatorname{vec}\mathbf{B})$ and $\operatorname{vec}(\mathbf{ABC})=(\mathbf{C}'\otimes\mathbf{A})\operatorname{vec}\mathbf{B}$ yields

$$\sum_{i=1}^{n-1}(n-i)^{1/2}\operatorname{tr}[\mathbf{d}_i^{(n)\prime}(\boldsymbol{\theta}_0)\boldsymbol{\Gamma}^{(n)}_{i;\boldsymbol{\Sigma},f}(\boldsymbol{\theta}_0)]=\begin{pmatrix}\mathbf{a}_1^{(n)}+\mathbf{b}_1^{(n)}\\ \vdots\\ \mathbf{a}_{n-1}^{(n)}+\mathbf{b}_{n-1}^{(n)}\end{pmatrix}'\mathbf{S}^{(n)}_{\boldsymbol{\Sigma},f}(\boldsymbol{\theta}_0), \tag{26}$$

with $\mathbf{a}_i^{(n)}$ and $\mathbf{b}_i^{(n)}$ defined at the end of Section 5.

Since $\mathbf{H}_m^{\prime(l)}\mathbf{B}_m^{(r)}\mathbf{H}_{m,q_1}^{(r)}=\mathbf{H}_{m,q_1}^{\prime(l)}$, (26) can be written as

$$\begin{aligned}\boldsymbol{\tau}^{(n)\prime}\boldsymbol{\Delta}^{(n)}_{\boldsymbol{\Sigma},f}(\boldsymbol{\theta}_0)&=[(\mathbf{H}^{(r)}_{n-1}\mathbf{B}^{\prime(l)}_{n-1}\mathbf{G}^{\prime(l)}_{n-1,p_1}|\mathbf{H}^{(r)}_{n-1,q_1})\boldsymbol{\tau}^{(n)}]'\mathbf{S}^{(n)}_{\boldsymbol{\Sigma},f}(\boldsymbol{\theta}_0)\\&=[(\mathbf{G}^{\prime(l)}_{n-1,p_1}|\mathbf{H}^{\prime(l)}_{n-1,q_1})\boldsymbol{\tau}^{(n)}]'(\mathbf{H}^{(r)}_{n-1}\mathbf{B}^{\prime(l)}_{n-1})'\mathbf{S}^{(n)}_{\boldsymbol{\Sigma},f}(\boldsymbol{\theta}_0)\\&=\begin{pmatrix}\tilde{\mathbf{a}}_1^{(n)}+\tilde{\mathbf{b}}_1^{(n)}\\ \vdots\\ \tilde{\mathbf{a}}_{n-1}^{(n)}+\tilde{\mathbf{b}}_{n-1}^{(n)}\end{pmatrix}'(\mathbf{H}^{(r)}_{n-1}\mathbf{B}^{\prime(l)}_{n-1})'\mathbf{S}^{(n)}_{\boldsymbol{\Sigma},f}(\boldsymbol{\theta}_0),\end{aligned} \tag{27}$$

where $\tilde{\mathbf{a}}_t^{(n)}:=\sum_{j=1}^{\min(p_1,t)}(\mathbf{G}'_{t-j}\otimes\mathbf{I}_k)(\operatorname{vec}\boldsymbol{\gamma}_j^{(n)})$ and $\tilde{\mathbf{b}}_t^{(n)}:=\sum_{j=1}^{\min(q_1,t)}(\mathbf{H}'_{t-j}\otimes\mathbf{I}_k)(\operatorname{vec}\boldsymbol{\delta}_j^{(n)})$; the sequences $(\tilde{\mathbf{a}}_t^{(n)})$ and $(\tilde{\mathbf{b}}_t^{(n)})$ clearly satisfy

$$(\tilde{\mathbf{a}}_1^{(n)\prime}+\tilde{\mathbf{b}}_1^{(n)\prime},\dots,\tilde{\mathbf{a}}_{\pi_0}^{(n)\prime}+\tilde{\mathbf{b}}_{\pi_0}^{(n)\prime})'=\mathbf{M}_{\boldsymbol{\theta}_0}\boldsymbol{\tau}^{(n)}. \tag{28}$$

Note that, for $t\ge p_0+q_0+1$,

$$\mathbf{D}(L)\mathbf{G}'_t=\mathbf{D}(L)\left(\sum_{i=1}^{p_0}\mathbf{G}'_{t-i}\mathbf{A}'_i\right)=\sum_{i=1}^{p_0}(\mathbf{D}(L)\mathbf{G}'_{t-i})\mathbf{A}'_i=\mathbf{0}.$$

Therefore, $\mathbf{D}(L)\mathbf{G}'_t=\mathbf{0}$ for all $t\ge q_0+1$. In the same way, we obtain that $\mathbf{D}(L)\mathbf{H}'_t=\mathbf{0}$ for $t\ge p_0+1$. Now, consider the $k^2$-dimensional operator $\mathbf{D}^{(l)}(L):=\mathbf{I}_{k^2}+\sum_{i=1}^{p_0+q_0}(\mathbf{D}_i\otimes\mathbf{I}_k)L^i$. This operator is such that, for $t-p_1\ge q_0+1$, $\mathbf{D}^{(l)}(L)\tilde{\mathbf{a}}_t^{(n)} = \sum_{j=1}^{p_1}(\mathbf{D}(L)\mathbf{G}'_{t-j} \otimes \mathbf{I}_k)(\operatorname{vec}\boldsymbol{\gamma}_j^{(n)}) = \mathbf{0}$. Similarly, one can check that for $t - q_1 \geq p_0 + 1$, $\mathbf{D}^{(l)}(L)\tilde{\mathbf{b}}_t^{(n)} = \sum_{j=1}^{q_1}(\mathbf{D}(L)\mathbf{H}'_{t-j} \otimes \mathbf{I}_k)(\operatorname{vec}\boldsymbol{\delta}_j^{(n)}) = \mathbf{0}$. This implies that $\tilde{\mathbf{a}}_t^{(n)} + \tilde{\mathbf{b}}_t^{(n)}$ satisfies $\mathbf{D}^{(l)}(L)(\tilde{\mathbf{a}}_t^{(n)} + \tilde{\mathbf{b}}_t^{(n)}) = \mathbf{0}$ for all $t \geq \max(p_1 + q_0 + 1, q_1 + p_0 + 1) = \pi + (p_0 + q_0) + 1$. Since $\{\boldsymbol{\Psi}_t^{(1)} \otimes \mathbf{I}_k, \ldots, \boldsymbol{\Psi}_t^{(p_0+q_0)} \otimes \mathbf{I}_k\}$ is a fundamental system of solutions of the homogeneous difference equation associated with $\mathbf{D}^{(l)}(L)$, we have

$$\begin{pmatrix} \tilde{\mathbf{a}}_{\pi+1}^{(n)} + \tilde{\mathbf{b}}_{\pi+1}^{(n)} \\ \vdots \\ \tilde{\mathbf{a}}_{n-1}^{(n)} + \tilde{\mathbf{b}}_{n-1}^{(n)} \end{pmatrix} = \bar{\boldsymbol{\Psi}}_{n-1}\mathbf{C}_{\Psi}^{-1} \begin{pmatrix} \tilde{\mathbf{a}}_{\pi+1}^{(n)} + \tilde{\mathbf{b}}_{\pi+1}^{(n)} \\ \vdots \\ \tilde{\mathbf{a}}_{\pi_0}^{(n)} + \tilde{\mathbf{b}}_{\pi_0}^{(n)} \end{pmatrix} \tag{29}$$

[see, e.g., Hallin (1986)]. Combining (28) and (29), we obtain

$$(\tilde{\mathbf{a}}_1^{(n)\prime} + \tilde{\mathbf{b}}_1^{(n)\prime}, \ldots, \tilde{\mathbf{a}}_{n-1}^{(n)\prime} + \tilde{\mathbf{b}}_{n-1}^{(n)\prime})' = \begin{pmatrix} \mathbf{I}_{k^2\pi} & \mathbf{0} \\ \mathbf{0} & \bar{\boldsymbol{\Psi}}_{n-1}\mathbf{C}_{\Psi}^{-1} \end{pmatrix} \mathbf{M}_{\boldsymbol{\theta}_0}\boldsymbol{\tau}^{(n)},$$

which, together with (27), establishes the result. □

PROOF OF LEMMA 2. Under $\bigcup_f \mathcal{H}^{(n)}(\boldsymbol{\theta}_0, \boldsymbol{\Sigma}, f)$, the residuals $\mathbf{Z}_1(\boldsymbol{\theta}_0), \ldots, \mathbf{Z}_n(\boldsymbol{\theta}_0)$, from which $\mathbf{C}_{\mathrm{Tyl}}^{(n)}$ is computed, are i.i.d. and elliptically symmetric, with mean $\mathbf{0}$ and scatter matrix $\boldsymbol{\Sigma}$. Tyler (1987) showed that $\mathbf{C}_{\mathrm{Tyl}}^{(n)}$ then is root-$n$ consistent for $\mathbf{C}_0 := c^{-1}\boldsymbol{\Sigma}^{-1/2}$, where $c$ denotes the upper left element in $\boldsymbol{\Sigma}^{-1/2}$. The result follows, since for any random vector $\mathbf{X}$,

$$\begin{aligned} &\left\| \frac{\mathbf{C}_{\mathrm{Tyl}}^{(n)}\mathbf{X}}{\|\mathbf{C}_{\mathrm{Tyl}}^{(n)}\mathbf{X}\|} - \frac{\boldsymbol{\Sigma}^{-1/2}\mathbf{X}}{\|\boldsymbol{\Sigma}^{-1/2}\mathbf{X}\|} \right\| \\ &\quad \leq \left| \frac{1}{\|\mathbf{C}_{\mathrm{Tyl}}^{(n)}\mathbf{X}\|} - \frac{1}{\|\mathbf{C}_0\mathbf{X}\|} \right| \|\mathbf{C}_{\mathrm{Tyl}}^{(n)}\mathbf{X}\| + \frac{1}{\|\mathbf{C}_0\mathbf{X}\|}\|\mathbf{C}_{\mathrm{Tyl}}^{(n)}\mathbf{X} - \mathbf{C}_0\mathbf{X}\| \\ &\quad \leq 2\frac{\|\mathbf{C}_{\mathrm{Tyl}}^{(n)}\mathbf{X} - \mathbf{C}_0\mathbf{X}\|}{\|\mathbf{C}_0\mathbf{X}\|} \leq 2\frac{\|\mathbf{C}_{\mathrm{Tyl}}^{(n)} - \mathbf{C}_0\|_{\mathcal{L}}\|\mathbf{X}\|}{\|\mathbf{C}_0\mathbf{X}\|} \leq 2\|\mathbf{C}_{\mathrm{Tyl}}^{(n)} - \mathbf{C}_0\|_{\mathcal{L}}\|\mathbf{C}_0^{-1}\|_{\mathcal{L}}, \end{aligned}$$

where $\|\mathbf{T}\|_{\mathcal{L}} := \sup\{\|\mathbf{T}\mathbf{x}\| \,|\, \|\mathbf{x}\| = 1\}$ denotes the operator norm of the square matrix $\mathbf{T}$. □

PROOF OF LEMMA 3. The definition of $\mathbf{W}_t^{(n)}(\mathbf{M})$ and (12) directly yield

$$\mathbf{W}_t^{(n)}(\mathbf{M}) = \frac{\mathbf{C}_{\mathrm{Tyl}}^{(n)}(\mathbf{M})\mathbf{M}\mathbf{Z}_t^{(n)}}{\|\mathbf{C}_{\mathrm{Tyl}}^{(n)}(\mathbf{M})\mathbf{M}\mathbf{Z}_t^{(n)}\|} = \frac{d\mathbf{O}\mathbf{C}_{\mathrm{Tyl}}^{(n)}\mathbf{Z}_t^{(n)}}{\|d\mathbf{O}\mathbf{C}_{\mathrm{Tyl}}^{(n)}\mathbf{Z}_t^{(n)}\|} = \mathbf{O}\mathbf{W}_t^{(n)}.$$

□

**A.2. Proofs of Propositions 2 and 3.** The following lemma, which follows along the same lines as Lemma 4 in Hallin and Paindaveine (2002b), will be used in the proof of Proposition 2.

LEMMA 5. *Let $i \in \{1, \ldots, n-1\}$ and $t, \tilde{t} \in \{i+1, \ldots, n\}$ be such that $t \neq \tilde{t}$. Assume that $g : \mathbb{R}^{nk} = \mathbb{R}^k \times \cdots \times \mathbb{R}^k \to \mathbb{R}$ is even in all its arguments, and such that the expectation below exists. Then, under $\bigcup_f \mathcal{H}^{(n)}(\boldsymbol{\theta}_0, \boldsymbol{\Sigma}, f)$,*

$$\mathrm{E}[g(\mathbf{Z}_1^{(n)}(\boldsymbol{\theta}_0), \ldots, \mathbf{Z}_n^{(n)}(\boldsymbol{\theta}_0))(\mathbf{P}_t'\mathbf{Q}_{\tilde{t}})(\mathbf{R}_{t-i}'\mathbf{S}_{\tilde{t}-i})] = 0,$$

*where $\mathbf{P}_j, \mathbf{Q}_j, \mathbf{R}_j$ and $\mathbf{S}_j$ are any four statistics among $\mathbf{W}_j^{(n)}(\boldsymbol{\theta}_0)$ and $\mathbf{W}_j^{(n)}(\boldsymbol{\theta}_0) - \mathbf{U}_j^{(n)}(\boldsymbol{\theta}_0, \boldsymbol{\Sigma})$.*

PROOF OF PROPOSITION 2. Throughout, we write $d_t^{(n)}$, $R_t^{(n)}$, $\hat{R}_t^{(n)}$, $\mathbf{W}_t^{(n)}$ and $\mathbf{U}_t^{(n)}$ for $d_t^{(n)}(\boldsymbol{\theta}_0, \boldsymbol{\Sigma})$, $R_t^{(n)}(\boldsymbol{\theta}_0, \boldsymbol{\Sigma})$, $\hat{R}_t^{(n)}(\boldsymbol{\theta}_0)$, $\mathbf{W}_t^{(n)}(\boldsymbol{\theta}_0)$ and $\mathbf{U}_t^{(n)}(\boldsymbol{\theta}_0, \boldsymbol{\Sigma})$, respectively; all convergences and mathematical expectations are taken as $n \to \infty$, under $\mathcal{H}^{(n)}(\boldsymbol{\theta}_0, \boldsymbol{\Sigma}, f)$. Decompose

$$(n-i)^{1/2}[(\mathbf{C}_{\mathrm{Tyl}}^{(n)} \otimes (\mathbf{C}_{\mathrm{Tyl}}^{(n)\prime})^{-1}) \operatorname{vec} \tilde{\boldsymbol{\Gamma}}_{i;K}^{(n)}(\boldsymbol{\theta}_0) - (\boldsymbol{\Sigma}^{-1/2} \otimes \boldsymbol{\Sigma}'^{1/2}) \operatorname{vec} \boldsymbol{\Gamma}_{i;K;\boldsymbol{\Sigma},f}^{(n)}(\boldsymbol{\theta}_0)]$$

into $\operatorname{vec}(\mathbf{T}_1^{(n)} + \mathbf{T}_2^{(n)} + \mathbf{T}_3^{(n)})$, where

$$\mathbf{T}_1^{(n)} := (n-i)^{-1/2} \sum_{t=i+1}^{n} \biggl( K_1\biggl(\frac{\hat{R}_t^{(n)}}{n+1}\biggr) K_2\biggl(\frac{\hat{R}_{t-i}^{(n)}}{n+1}\biggr) - K_1\biggl(\frac{R_t^{(n)}}{n+1}\biggr) K_2\biggl(\frac{R_{t-i}^{(n)}}{n+1}\biggr)\biggr) \mathbf{W}_t^{(n)} \mathbf{W}_{t-i}^{(n)\prime},$$

$$\mathbf{T}_2^{(n)} := (n-i)^{-1/2} \sum_{t=i+1}^{n} K_1\biggl(\frac{R_t^{(n)}}{n+1}\biggr) K_2\biggl(\frac{R_{t-i}^{(n)}}{n+1}\biggr) (\mathbf{W}_t^{(n)} \mathbf{W}_{t-i}^{(n)\prime} - \mathbf{U}_t^{(n)} \mathbf{U}_{t-i}^{(n)\prime})$$

and

$$\mathbf{T}_3^{(n)} := (n-i)^{-1/2} \sum_{t=i+1}^{n} \biggl( K_1\biggl(\frac{R_t^{(n)}}{n+1}\biggr) K_2\biggl(\frac{R_{t-i}^{(n)}}{n+1}\biggr) - K_1(\tilde{F}_k(d_t^{(n)})) K_2(\tilde{F}_k(d_{t-i}^{(n)}))\biggr) \mathbf{U}_t^{(n)} \mathbf{U}_{t-i}^{(n)\prime}.$$

We proceed by proving that $\operatorname{vec} \mathbf{T}_1^{(n)}, \operatorname{vec} \mathbf{T}_2^{(n)}$ and $\operatorname{vec} \mathbf{T}_3^{(n)}$ all converge to $\mathbf{0}$ in quadratic mean, as $n \to \infty$. Slutsky's classical argument then concludes the proof.

Let us start with $\mathbf{T}_3^{(n)}$. Using the fact that $(\operatorname{vec}\mathbf{A})'(\operatorname{vec}\mathbf{B}) = \operatorname{tr}(\mathbf{A}'\mathbf{B})$ and the independence between the $d_t$'s and the $\mathbf{U}_t$'s, we obtain

$$\|\operatorname{vec}\mathbf{T}_3^{(n)}\|_{L^2}^2 = \sum_{t=i+1}^{n} (c_{t;i}^{(n)})^2 \mathrm{E}\Biggl[\Biggl(K_1\biggl(\frac{R_t^{(n)}}{n+1}\biggr)K_2\biggl(\frac{R_{t-i}^{(n)}}{n+1}\biggr) - K_1(\tilde{F}_k(d_t^{(n)}))K_2(\tilde{F}_k(d_{t-i}^{(n)}))\Biggr)^2\Biggr],$$

where $c_{t;i}^{(n)} = (n-i)^{-1/2}$ for all $t = i+1, \dots, n$. Hájek's projection result thus implies that $\|\operatorname{vec}\mathbf{T}_3^{(n)}\|_{L^2}^2 = o(1)$ as $n \to \infty$. The same result also implies that, for all $t = i+1, \dots, n$,

$$\mathrm{E}\Biggl[\Biggl(K_1\biggl(\frac{R_t^{(n)}}{n+1}\biggr)K_2\biggl(\frac{R_{t-i}^{(n)}}{n+1}\biggr) - K_1(\tilde{F}_k(d_t^{(n)}))K_2(\tilde{F}_k(d_{t-i}^{(n)}))\Biggr)^2\Biggr] = o(1). \tag{30}$$

For $\mathbf{T}_2^{(n)}$, decomposing $\mathbf{W}_t^{(n)}\mathbf{W}_{t-i}^{(n)\prime} - \mathbf{U}_t^{(n)}\mathbf{U}_{t-i}^{(n)\prime}$ into $(\mathbf{W}_t^{(n)} - \mathbf{U}_t^{(n)})\mathbf{W}_{t-i}^{(n)\prime} + \mathbf{U}_t^{(n)}(\mathbf{W}_{t-i}^{(n)} - \mathbf{U}_{t-i}^{(n)})'$, then using the identity $(\operatorname{vec}\mathbf{A})'(\operatorname{vec}\mathbf{B}) = \operatorname{tr}(\mathbf{A}'\mathbf{B})$ again and Lemma 5, one obtains

$$\begin{aligned}
&\|\operatorname{vec}\mathbf{T}_2^{(n)}\|_{L^2}^2 \\
&\quad\le 2(n-i)^{-1}\sum_{t=i+1}^{n}\mathrm{E}\Biggl[\Biggl(K_1\biggl(\frac{R_t^{(n)}}{n+1}\biggr)K_2\biggl(\frac{R_{t-i}^{(n)}}{n+1}\biggr)\Biggr)^2\|\mathbf{W}_t^{(n)} - \mathbf{U}_t^{(n)}\|^2\Biggr] \\
&\qquad + 2(n-i)^{-1}\sum_{t=i+1}^{n}\mathrm{E}\Biggl[\Biggl(K_1\biggl(\frac{R_t^{(n)}}{n+1}\biggr)K_2\biggl(\frac{R_{t-i}^{(n)}}{n+1}\biggr)\Biggr)^2\|\mathbf{W}_{t-i}^{(n)} - \mathbf{U}_{t-i}^{(n)}\|^2\Biggr].
\end{aligned} \tag{31}$$

Consider the first term in the right-hand side of (31) (the second term can be dealt with in the same way). Let $A_{t;i}^{(n)} := K_1(R_t^{(n)}/(n+1))K_2(R_{t-i}^{(n)}/(n+1))$ and $B_{t;i}^{(n)} := K_1(\tilde{F}_k(d_t^{(n)}))K_2(\tilde{F}_k(d_{t-i}^{(n)}))$. Using (30) and the independence between the $d_t^{(n)}$'s and the $\mathbf{U}_t^{(n)}$'s, we obtain

$$\begin{aligned}
\mathrm{E}[(A_{t;i}^{(n)})^2\|\mathbf{W}_t^{(n)} - \mathbf{U}_t^{(n)}\|^2] &= \mathrm{E}[(B_{t;i}^{(n)})^2\|\mathbf{W}_t^{(n)} - \mathbf{U}_t^{(n)}\|^2] + o(1) \\
&= \|K_1(U)\|_{L^2}^2\,\|K_2(U)\|_{L^2}^2\|\mathbf{W}_t^{(n)} - \mathbf{U}_t^{(n)}\|_{L_2}^2 + o(1),
\end{aligned}$$

where $U$ is uniformly distributed over $]0,1[$. Lemma 2 thus implies that $\operatorname{vec}\mathbf{T}_2^{(n)}$ converges to $\mathbf{0}$ in quadratic mean.

Finally, using Lemma 5 again,

$$\|\operatorname{vec}\mathbf{T}_1^{(n)}\|_{L^2}^2 = (n-i)^{-1}\sum_{t=i+1}^{n}\Biggl(K_1\biggl(\frac{\hat{R}_t^{(n)}}{n+1}\biggr)K_2\biggl(\frac{\hat{R}_{t-i}^{(n)}}{n+1}\biggr) - K_1\biggl(\frac{R_t^{(n)}}{n+1}\biggr)K_2\biggl(\frac{R_{t-i}^{(n)}}{n+1}\biggr)\Biggr)^2.$$

This entails that $\operatorname{vec} \mathbf{T}_1^{(n)}$ also is $o_{\mathrm{qm}}(1)$, provided that

$$K_1\left(\frac{\hat{R}_t^{(n)}}{n+1}\right)K_2\left(\frac{\hat{R}_{t-i}^{(n)}}{n+1}\right) - K_1\left(\frac{R_t^{(n)}}{n+1}\right)K_2\left(\frac{R_{t-i}^{(n)}}{n+1}\right) \stackrel{L^2}{\to} 0 \qquad \text{as } n \to \infty. \tag{32}$$

Now, Lemma 1 establishes the same convergence as in (32), but in probability. On the other hand, it follows from (30) that $[K_1(R_t/(n+1))K_2(R_{t-i}/(n+1))]^2$ is uniformly integrable, which (in view of the invariance of Tyler's estimator of scatter under permutations of the residuals) implies that $[K_1(\hat{R}_t/(n+1))K_2(\hat{R}_{t-i}/(n+1))]^2$ also is. The $L^2$ convergence in (32) follows.

Summing up, we have shown that $(n-i)^{1/2}[(\mathbf{C}_{\mathrm{Tyl}}^{(n)} \otimes (\mathbf{C}_{\mathrm{Tyl}}^{(n)\prime})^{-1}) \operatorname{vec} \widetilde{\mathbf{\Gamma}}_{i;K}^{(n)}(\boldsymbol{\theta}_0) - (\mathbf{\Sigma}^{-1/2} \otimes \mathbf{\Sigma}'^{1/2}) \operatorname{vec} \mathbf{\Gamma}_{i;K;\mathbf{\Sigma},f}^{(n)}(\boldsymbol{\theta}_0)]$ is $o_{\mathrm{qm}}(1)$ as $n \to \infty$. This concludes the proof, since, from a multivariate application of Slutsky's theorem,

$$\begin{aligned}(n-i)^{1/2}[(\mathbf{C}_{\mathrm{Tyl}}^{(n)} \otimes (\mathbf{C}_{\mathrm{Tyl}}^{(n)\prime}))^{-1} \operatorname{vec} \widetilde{\mathbf{\Gamma}}_{i;K}^{(n)}(\boldsymbol{\theta}_0) \\ - (\mathbf{\Sigma}^{-1/2} \otimes \mathbf{\Sigma}'^{1/2}) \operatorname{vec} \widetilde{\mathbf{\Gamma}}_{i;K}^{(n)}(\boldsymbol{\theta}_0)] = o_{\mathrm{P}}(1),\end{aligned}$$

under $\mathcal{H}^{(n)}(\boldsymbol{\theta}_0, \mathbf{\Sigma}, f)$. □

Proof of Proposition 3. Under $\mathcal{H}^{(n)}(\boldsymbol{\theta}_0, \mathbf{\Sigma}, f)$, one can use the same argument as in Lemma 4.12 in Garel and Hallin (1995). The result under the sequence of alternatives is obtained as usual, first establishing the joint normality of $\mathbf{S}_{m;K,\mathbf{\Sigma},f}^{(n)}(\boldsymbol{\theta}_0)$ and $L_{\boldsymbol{\theta}_0+n^{-1/2}\boldsymbol{\tau}/\boldsymbol{\theta}_0;\mathbf{\Sigma},f}^{(n)}$ under $\mathcal{H}^{(n)}(\boldsymbol{\theta}_0, \mathbf{\Sigma}, f)$, then applying Le Cam's third Lemma; the required joint normality easily follows from a routine application of the classical Cramér–Wold device. □

### A.3. Proofs of Propositions 4 and 5.

Proof of Proposition 4. (i) Let $\{\mathbf{\Psi}_t^{(1)}, \ldots, \mathbf{\Psi}_t^{(p_0+q_0)}\}$ and $\{\mathbf{\Phi}_t^{(1)}, \ldots, \mathbf{\Phi}_t^{(p_0+q_0)}\}$ be two fundamental systems of solutions associated with $\mathbf{D}(L)$. The vector structure of the space of solutions of $\mathbf{D}(L)\mathbf{x}_t = \mathbf{0}$, $\mathbf{x}_t \in \mathbb{R}^k$ implies that, for all $j = 1, \ldots, p_0 + q_0$, there exists a $k(p_0+q_0) \times k$ matrix $\mathbf{\Lambda}_j$ such that $\mathbf{\Phi}_t^{(j)} = (\mathbf{\Psi}_t^{(1)}, \ldots, \mathbf{\Psi}_t^{(p_0+q_0)})\mathbf{\Lambda}_j$. Letting $\mathbf{\Lambda} := (\mathbf{\Lambda}_1, \ldots, \mathbf{\Lambda}_{p_0+q_0})$, this implies that

$$\begin{pmatrix} \mathbf{\Phi}_{\pi+1}^{(1)} & \cdots & \mathbf{\Phi}_{\pi+1}^{(p_0+q_0)} \\ \mathbf{\Phi}_{\pi+2}^{(1)} & \cdots & \mathbf{\Phi}_{\pi+2}^{(p_0+q_0)} \\ \vdots & & \vdots \\ \mathbf{\Phi}_{m}^{(1)} & \cdots & \mathbf{\Phi}_{m}^{(p_0+q_0)} \end{pmatrix} = \begin{pmatrix} \mathbf{\Psi}_{\pi+1}^{(1)} & \cdots & \mathbf{\Psi}_{\pi+1}^{(p_0+q_0)} \\ \mathbf{\Psi}_{\pi+2}^{(1)} & \cdots & \mathbf{\Psi}_{\pi+2}^{(p_0+q_0)} \\ \vdots & & \vdots \\ \mathbf{\Psi}_{m}^{(1)} & \cdots & \mathbf{\Psi}_{m}^{(p_0+q_0)} \end{pmatrix} \mathbf{\Lambda},$$

so that $\bar{\mathbf{\Phi}}_m = \bar{\mathbf{\Psi}}_m(\mathbf{\Lambda} \otimes \mathbf{I}_k)$, where $\bar{\mathbf{\Phi}}_m$ is the equivalent of $\bar{\mathbf{\Psi}}_m$, but computed from the $\mathbf{\Phi}_t^{(j)}$'s. Thus, with obvious notation, $\mathbf{Q}_{\boldsymbol{\theta}_0;\mathbf{\Phi}}^{(m)} = \mathbf{Q}_{\boldsymbol{\theta}_0;\mathbf{\Psi}}^{(m)}\bar{\mathbf{\Lambda}}$ for all $m$, where $\bar{\mathbf{\Lambda}} := \binom{\mathbf{I}_{k^2\pi} \quad \mathbf{0}}{\mathbf{0} \quad \mathbf{\Lambda}\otimes\mathbf{I}_k}$, yielding (note that since the $\mathbf{\Psi}_t^{(j)}$'s and $\mathbf{\Phi}_t^{(j)}$'s constitute fundamental systems, $\mathbf{\Lambda}$, and hence $\bar{\mathbf{\Lambda}}$, are nonsingular)

$$\begin{aligned}\widetilde{\mathbf{T}}^{(n)\prime}_{K;\mathbf{\Phi}}(\boldsymbol{\theta}_0)(\mathbf{J}^{(n)}_{\boldsymbol{\theta}_0,\widehat{\mathbf{\Sigma}};\mathbf{\Phi}})^{-1}\widetilde{\mathbf{T}}^{(n)}_{K;\mathbf{\Phi}}(\boldsymbol{\theta}_0) &= [\widetilde{\mathbf{T}}^{(n)\prime}_{K;\mathbf{\Psi}}(\boldsymbol{\theta}_0)\bar{\mathbf{\Lambda}}][\bar{\mathbf{\Lambda}}'\mathbf{J}^{(n)}_{\boldsymbol{\theta}_0,\widehat{\mathbf{\Sigma}};\mathbf{\Psi}}\bar{\mathbf{\Lambda}}]^{-1}[\bar{\mathbf{\Lambda}}'\widetilde{\mathbf{T}}^{(n)}_{K;\mathbf{\Psi}}(\boldsymbol{\theta}_0)]\\ &= \widetilde{\mathbf{T}}^{(n)\prime}_{K;\mathbf{\Psi}}(\boldsymbol{\theta}_0)(\mathbf{J}^{(n)}_{\boldsymbol{\theta}_0,\widehat{\mathbf{\Sigma}};\mathbf{\Psi}})^{-1}\widetilde{\mathbf{T}}^{(n)}_{K;\mathbf{\Psi}}(\boldsymbol{\theta}_0),\end{aligned}$$

as was to be proved. The statement about the dependence on $p_1$ and $q_1$ is trivial, since $\widetilde{\mathbf{T}}^{(n)}_K(\boldsymbol{\theta}_0)$, $\mathbf{J}^{(n)}_{\boldsymbol{\theta}_0,\widehat{\mathbf{\Sigma}}}$, as well as $\pi_0$, depend on $p_1$ and $q_1$ only through $\pi$.

(ii) Letting

$$n^{1/2}\widetilde{\mathbf{T}}^{(n)}_{K;\mathbf{\Sigma}}(\boldsymbol{\theta}_0) := \mathbf{Q}^{(n)\prime}_{\boldsymbol{\theta}_0}((n-1)^{1/2}(\operatorname{vec}\widetilde{\mathbf{\Gamma}}^{(n)}_{1;K;\mathbf{\Sigma}}(\boldsymbol{\theta}_0))', \ldots, (\operatorname{vec}\widetilde{\mathbf{\Gamma}}^{(n)}_{n-1;K;\mathbf{\Sigma}}(\boldsymbol{\theta}_0))')',$$

with

$$\begin{aligned}\widetilde{\mathbf{\Gamma}}^{(n)}_{i;K;\mathbf{\Sigma}}(\boldsymbol{\theta}_0) := \mathbf{C}^{(n)\prime}_{\mathrm{Tyl}}\Bigg(\frac{1}{n-i}\sum_{t=i+1}^{n} K_1\bigg(\frac{R^{(n)}_t(\boldsymbol{\theta}_0,\mathbf{\Sigma})}{n+1}\bigg)K_2\bigg(\frac{R^{(n)}_{t-i}(\boldsymbol{\theta}_0,\mathbf{\Sigma})}{n+1}\bigg)\\ \times \mathbf{W}^{(n)}_t(\boldsymbol{\theta}_0)\mathbf{W}^{(n)\prime}_{t-i}(\boldsymbol{\theta}_0)\Bigg)(\mathbf{C}^{(n)\prime}_{\mathrm{Tyl}})^{-1},\end{aligned}$$

one can verify (proceeding as for the first term in the decomposition argument in the proof of Proposition 2) that $n^{1/2}(\widetilde{\mathbf{T}}^{(n)}_K(\boldsymbol{\theta}_0) - \widetilde{\mathbf{T}}^{(n)}_{K;\mathbf{\Sigma}}(\boldsymbol{\theta}_0))$ tends to zero in quadratic mean as $n\to\infty$ under $\bigcup_f \mathcal{H}^{(n)}(\boldsymbol{\theta}_0,\mathbf{\Sigma},f)$. This entails that

$$\begin{aligned}Q^{(n)}_K(\boldsymbol{\theta}_0) &= \frac{k^2}{\mathrm{E}[K_1^2(U)]\mathrm{E}[K_2^2(U)]}\\ &\quad\times (n^{1/2}\widetilde{\mathbf{T}}^{(n)}_{K;\mathbf{\Sigma}}(\boldsymbol{\theta}_0))'(\mathbf{J}^{(n)}_{\boldsymbol{\theta}_0,\widehat{\mathbf{\Sigma}}})^{-1}(n^{1/2}\widetilde{\mathbf{T}}^{(n)}_{K;\mathbf{\Sigma}}(\boldsymbol{\theta}_0)) + o_{\mathrm{P}}(1)\end{aligned}$$

is asymptotically invariant with respect to $\mathcal{G}^{(n)}_{\mathbf{\Sigma}}$ under $\bigcup_f \mathcal{H}^{(n)}(\boldsymbol{\theta}_0,\mathbf{\Sigma},f)$, since $n^{1/2}\widetilde{\mathbf{T}}^{(n)}_{K;\mathbf{\Sigma}}(\boldsymbol{\theta}_0)$ and $\mathbf{J}^{(n)}_{\boldsymbol{\theta}_0,\widehat{\mathbf{\Sigma}}}$ are strictly invariant with respect to the same group.

(iii), (iv) Proposition 2 and the multivariate Slutsky theorem show that $Q^{(n)}_K(\boldsymbol{\theta}_0)$ has the same asymptotic behavior [under $\mathcal{H}^{(n)}(\boldsymbol{\theta}_0,\mathbf{\Sigma},f)$, as well as under the sequence of local alternatives $\mathcal{H}^{(n)}(\boldsymbol{\theta}_0+n^{-1/2}\boldsymbol{\tau},\mathbf{\Sigma},f)$] as

$$\frac{k^2}{\mathrm{E}[K_1^2(U)]\mathrm{E}[K_2^2(U)]}(n^{1/2}\mathbf{T}^{(n)}_{K;\mathbf{\Sigma},f}(\boldsymbol{\theta}_0))'\mathbf{J}^{-1}_{\boldsymbol{\theta}_0,\mathbf{\Sigma}}(n^{1/2}\mathbf{T}^{(n)}_{K;\mathbf{\Sigma},f}(\boldsymbol{\theta}_0)),$$

where $n^{1/2}\mathbf{T}^{(n)}_{K;\mathbf{\Sigma},f}(\boldsymbol{\theta}_0) := \mathbf{Q}^{(n)\prime}_{\boldsymbol{\theta}_0}\mathbf{S}^{(n)}_{n-1;K;\mathbf{\Sigma},f}(\boldsymbol{\theta}_0)$ [see (16)]. Now, Proposition 3 and a classical result on triangular arrays [Brockwell and Davis (1987),

Proposition 6.3.9] imply that $n^{1/2}\mathbf{T}^{(n)}_{K;\boldsymbol{\Sigma},f}(\boldsymbol{\theta}_0)$ is asymptotically $k^2\pi_0$-variate normal, with mean $\mathbf{0}$ under $\mathcal{H}^{(n)}(\boldsymbol{\theta}_0, \boldsymbol{\Sigma}, f)$, and mean

$$\frac{1}{k^2} D_k(K_2; f) C_k(K_1; f) \mathbf{J}_{\boldsymbol{\theta}_0,\boldsymbol{\Sigma}} \mathbf{P}_{\boldsymbol{\theta}_0} \mathbf{M}_{\boldsymbol{\theta}_0} \boldsymbol{\tau}$$

under $\mathcal{H}^{(n)}(\boldsymbol{\theta}_0 + n^{-1/2}\boldsymbol{\tau}, \boldsymbol{\Sigma}, f)$, and with covariance matrix $(\mathrm{E}[K_1^2(U)] \times \mathrm{E}[K_2^2(U)]/k^2)\mathbf{J}_{\boldsymbol{\theta}_0,\boldsymbol{\Sigma}}$ under both. The result follows.

(v) It follows from Le Cam [(1986), Section 11.9] and the LAN property in Proposition 1 that the test $\underline{\phi}^{(n)}_{\boldsymbol{\Sigma},f_\star}$ rejecting the null hypothesis whenever

$$\boldsymbol{\Delta}^{(n)\prime}_{\boldsymbol{\Sigma},f_\star}(\boldsymbol{\theta}_0)(\boldsymbol{\Gamma}_{\boldsymbol{\Sigma},f_\star}(\boldsymbol{\theta}_0))^{-}\boldsymbol{\Delta}^{(n)}_{\boldsymbol{\Sigma},f_\star}(\boldsymbol{\theta}_0) > \chi^2_{s,1-\alpha},$$

where $\mathbf{A}^-$ denotes any arbitrary generalized inverse of $\mathbf{A}$ and $s := \mathrm{rank}(\boldsymbol{\Gamma}_{\boldsymbol{\Sigma},f_\star}(\boldsymbol{\theta}_0))$, is locally and asymptotically maximin, at probability level $\alpha$, for $\mathcal{H}^{(n)}(\boldsymbol{\theta}_0, \boldsymbol{\Sigma}, f_\star)$ against $\bigcup_{\boldsymbol{\theta}\neq\boldsymbol{\theta}_0} \mathcal{H}^{(n)}(\boldsymbol{\theta}, \boldsymbol{\Sigma}, f_\star)$. Note that $\mathrm{rank}(\boldsymbol{\Gamma}_{\boldsymbol{\Sigma},f_\star}(\boldsymbol{\theta}_0)) = \mathrm{rank}(\mathbf{M}'_{\boldsymbol{\theta}_0}\mathbf{P}'_{\boldsymbol{\theta}_0} \times \mathbf{J}_{\boldsymbol{\theta}_0,\boldsymbol{\Sigma}}\mathbf{P}_{\boldsymbol{\theta}_0}\mathbf{M}_{\boldsymbol{\theta}_0}) = \min(k^2(p_1+q_1), k^2\pi_0) = k^2\pi_0$, since $\mathbf{M}_{\boldsymbol{\theta}_0}$, $\mathbf{P}_{\boldsymbol{\theta}_0}$ and $\mathbf{J}_{\boldsymbol{\theta}_0,\boldsymbol{\Sigma}}$ have maximal rank. Of course, the same optimality property holds for the asymptotically equivalent [under $\mathcal{H}^{(n)}(\boldsymbol{\theta}_0, \boldsymbol{\Sigma}, f_\star)$, as well as under contiguous alternatives] test $\underline{\phi}^{(n)}_{f_\star}$ that rejects the null hypothesis whenever

$$\boldsymbol{\Delta}^{(n)\prime}_{K_{f_\star}}(\boldsymbol{\theta}_0)(\hat{\boldsymbol{\Gamma}}^{(n)}_{f_\star}(\boldsymbol{\theta}_0))^{-}\boldsymbol{\Delta}^{(n)}_{K_{f_\star}}(\boldsymbol{\theta}_0) > \chi^2_{k^2\pi_0,1-\alpha},$$

where $\boldsymbol{\Delta}^{(n)}_{K_{f_\star}}(\boldsymbol{\theta}_0) := n^{1/2}\mathbf{M}'_{\boldsymbol{\theta}_0}\mathbf{P}'_{\boldsymbol{\theta}_0}\widetilde{\mathbf{T}}^{(n)}_K(\boldsymbol{\theta}_0)$, with $K_1 := \varphi_{f_\star}\circ\tilde{F}^{-1}_{\star k}$ and $K_2 = \tilde{F}^{-1}_{\star k}$, and

$$\hat{\boldsymbol{\Gamma}}^{(n)}_{f_\star}(\boldsymbol{\theta}_0) := \frac{\mu_{k+1;f_\star}\mathcal{I}_{k,f_\star}}{k^2\mu_{k-1;f_\star}}\mathbf{M}'_{\boldsymbol{\theta}_0}\mathbf{P}'_{\boldsymbol{\theta}_0}\mathbf{J}^{(n)}_{\boldsymbol{\theta}_0,\widehat{\boldsymbol{\Sigma}}}\mathbf{P}_{\boldsymbol{\theta}_0}\mathbf{M}_{\boldsymbol{\theta}_0} = \boldsymbol{\Gamma}_{\boldsymbol{\Sigma},f_\star}(\boldsymbol{\theta}_0) + o_{\mathrm{P}}(1)$$

under $\mathcal{H}^{(n)}(\boldsymbol{\theta}_0, \boldsymbol{\Sigma}, f_\star)$. But, in view of Lemma 2.2.5(c) of Rao and Mitra (1971),

$$\begin{aligned}
&\boldsymbol{\Delta}^{(n)\prime}_{K_{f_\star}}(\boldsymbol{\theta}_0)(\hat{\boldsymbol{\Gamma}}^{(n)}_{f_\star}(\boldsymbol{\theta}_0))^{-}\boldsymbol{\Delta}^{(n)}_{K_{f_\star}}(\boldsymbol{\theta}_0)\\
&\quad= \frac{k^2 n\mu_{k-1;f_\star}}{\mu_{k+1;f_\star}\mathcal{I}_{k,f_\star}}\widetilde{\mathbf{T}}^{(n)\prime}_K(\boldsymbol{\theta}_0)\mathbf{P}_{\boldsymbol{\theta}_0}\mathbf{M}_{\boldsymbol{\theta}_0}(\mathbf{M}'_{\boldsymbol{\theta}_0}\mathbf{P}'_{\boldsymbol{\theta}_0}\mathbf{J}^{(n)}_{\boldsymbol{\theta}_0,\widehat{\boldsymbol{\Sigma}}}\mathbf{P}_{\boldsymbol{\theta}_0}\mathbf{M}_{\boldsymbol{\theta}_0})^{-}\mathbf{M}'_{\boldsymbol{\theta}_0}\mathbf{P}'_{\boldsymbol{\theta}_0}\widetilde{\mathbf{T}}^{(n)}_K(\boldsymbol{\theta}_0)\\
&\quad= \frac{k^2 n}{\mathrm{E}[(\varphi_{f_\star}(\tilde{F}^{-1}_{\star k}(U)))^2]\mathrm{E}[(\tilde{F}^{-1}_{\star k}(U))^2]}\widetilde{\mathbf{T}}^{(n)\prime}_K(\boldsymbol{\theta}_0)(\mathbf{J}^{(n)}_{\boldsymbol{\theta}_0,\widehat{\boldsymbol{\Sigma}}})^{-1}\widetilde{\mathbf{T}}^{(n)}_K(\boldsymbol{\theta}_0);
\end{aligned}$$

$\underline{\phi}^{(n)}_{f_\star}$ and $\phi^{(n)}_{f_\star}$, thus, are the same test. The result follows. □

PROOF OF PROPOSITION 5. (i) Model (1) under $\mathcal{H}^{(n)}(\boldsymbol{\theta}_0)$ can be written in the form

$$\mathbf{M}\mathbf{A}(L)\mathbf{M}^{-1}\mathbf{M}\mathbf{X}_t = \mathbf{M}\mathbf{B}(L)\mathbf{M}^{-1}\mathbf{M}\boldsymbol{\varepsilon}_t,$$

where $\mathbf{M}$ is an arbitrary full-rank $k \times k$ matrix. This null hypothesis is thus invariant under the group of affine transformations $\boldsymbol{\varepsilon}_t \mapsto \mathbf{M}\boldsymbol{\varepsilon}_t$ if and only if $\mathbf{M}\mathbf{A}_i\mathbf{M}^{-1} = \mathbf{A}_i$ for all $i = 1, \ldots, p_0$ and $\mathbf{M}\mathbf{B}_j\mathbf{M}^{-1} = \mathbf{B}_j$ for all $j = 1, \ldots, q_0$, that is, iff each $\mathbf{A}_i$ and each $\mathbf{B}_j$ commutes with any invertible matrix $\mathbf{M}$, which holds true iff they are proportional to the $k \times k$ identity matrix.

(ii) Let $\mathbf{M}$ be some nonsingular $k \times k$ matrix. For any statistic $T = T(\mathbf{X}^{(n)}_{-p_0+1}, \ldots, \mathbf{X}^{(n)}_n)$, write $T(\mathbf{M}) := T(\mathbf{M}\mathbf{X}^{(n)}_{-p_0+1}, \ldots, \mathbf{M}\mathbf{X}^{(n)}_n)$. It follows from Lemma 3 and from the equivariance properties of $\mathbf{C}^{(n)}_{\mathrm{Tyl}}$ that $\widetilde{\boldsymbol{\Gamma}}^{(n)}_{i;K}(\mathbf{M}) = \mathbf{M}'^{-1}\widetilde{\boldsymbol{\Gamma}}^{(n)}_{i;K}\mathbf{M}'$. Hence, $\widetilde{\mathbf{S}}^{(n)}_K(\mathbf{M}) = [\mathbf{I}_{n-1} \otimes (\mathbf{M} \otimes \mathbf{M}'^{-1})]\widetilde{\mathbf{S}}^{(n)}_K$. In the same way,

$$[\mathbf{I}_{n-1} \otimes (\widehat{\boldsymbol{\Sigma}}^{(n)}(\mathbf{M}) \otimes (\widehat{\boldsymbol{\Sigma}}^{(n)}(\mathbf{M}))^{-1})] = [\mathbf{I}_{n-1} \otimes (\mathbf{M} \otimes \mathbf{M}'^{-1})][\mathbf{I}_{n-1} \otimes (\widehat{\boldsymbol{\Sigma}}^{(n)} \otimes (\widehat{\boldsymbol{\Sigma}}^{(n)})^{-1})][\mathbf{I}_{n-1} \otimes (\mathbf{M} \otimes \mathbf{M}'^{-1})]'.$$

Now, $\mathbf{A}_i = a_i\mathbf{I}_k$ clearly implies that the Green matrices of the operator $\mathbf{A}(L)$ all are proportional to the identity matrix. The same property holds for $\mathbf{B}(L)$. It is then easy to verify that the operator $\mathbf{D}(L)$ also is scalar (meaning that $\mathbf{D}_i$ is proportional to the identity matrix for all $i = 1, \ldots, p_0 + q_0$). This implies that the fundamental system of solutions provided by Green's matrices of $\mathbf{D}(L)$ contains only matrices that are proportional to the identity matrix. Hence, $\mathbf{Q}^{(n)}_{\boldsymbol{\theta}_0} = \mathbf{W}^{(n)} \otimes \mathbf{I}_{k^2}$ for some $(n-1) \times \pi_0$ matrix $\mathbf{W}^{(n)}$. It follows that

$$[\mathbf{I}_{n-1} \otimes (\mathbf{M} \otimes \mathbf{M}'^{-1})]'\,\mathbf{Q}^{(n)}_{\boldsymbol{\theta}_0} = \mathbf{Q}^{(n)}_{\boldsymbol{\theta}_0}[\mathbf{I}_{\pi_0} \otimes (\mathbf{M} \otimes \mathbf{M}'^{-1})]',$$

which entails $\widetilde{\mathbf{T}}^{(n)}_K(\mathbf{M}) = [\mathbf{I}_{\pi_0} \otimes (\mathbf{M} \otimes \mathbf{M}'^{-1})]\widetilde{\mathbf{T}}^{(n)}_K$ and

$$\mathbf{J}^{(n)}_{\boldsymbol{\theta}_0,\widehat{\boldsymbol{\Sigma}}}(\mathbf{M}) = [\mathbf{I}_{\pi_0} \otimes (\mathbf{M} \otimes \mathbf{M}'^{-1})]\mathbf{J}^{(n)}_{\boldsymbol{\theta}_0,\widehat{\boldsymbol{\Sigma}}}[\mathbf{I}_{\pi_0} \otimes (\mathbf{M} \otimes \mathbf{M}'^{-1})]'.$$

Consequently, $Q^{(n)}_K(\mathbf{M}) = Q^{(n)}_K$. □

## REFERENCES

Bickel, P. J. (1982). On adaptive estimation. *Ann. Statist.* **10** 647–671. MR663424

Brockwell, P. J. and Davis, R. A. (1987). *Time Series: Theory and Methods*. Springer, New York. MR868859

Chernoff, H. and Savage, I. R. (1958). Asymptotic normality and efficiency of certain nonparametric test statistics. *Ann. Math. Statist.* **29** 972–994. MR100322

Drost, F. C., Klaassen, C. A. J. and Werker, B. J. M. (1997). Adaptive estimation in time-series models. *Ann. Statist.* **25** 786–818. MR1439324

Garel, B. and Hallin, M. (1995). Local asymptotic normality of multivariate ARMA processes with a linear trend. *Ann. Inst. Statist. Math.* **47** 551–579. MR1364260

Hájek, J. and Šidák, Z. (1967). *Theory of Rank Tests*. Academic Press, New York. MR229351

Hallin, M. (1986). Non-stationary $q$-dependent processes and time-varying moving-average models: Invertibility properties and the forecasting problem. *Adv. in Appl. Probab.* **18** 170–210. MR827335

Hallin, M. (1994). On the Pitman-nonadmissibility of correlogram-based methods. *J. Time Ser. Anal.* **15** 607–612. MR1312324

Hallin, M., Ingenbleek, J.-F. and Puri, M. L. (1985). Linear serial rank tests for randomness against ARMA alternatives. *Ann. Statist.* **13** 1156–1181. MR803764

Hallin, M., Ingenbleek, J.-F. and Puri, M. L. (1989). Asymptotically most powerful rank tests for multivariate randomness against serial dependence. *J. Multivariate Anal.* **30** 34–71. MR1003708

Hallin, M. and Mélard, G. (1988). Rank-based tests for randomness against first-order serial dependence. *J. Amer. Statist. Assoc.* **83** 1117–1128. MR997590

Hallin, M. and Paindaveine, D. (2002a). Optimal tests for multivariate location based on interdirections and pseudo-Mahalanobis ranks. *Ann. Statist.* **30** 1103–1133. MR1926170

Hallin, M. and Paindaveine, D. (2002b). Optimal procedures based on interdirections and pseudo-Mahalanobis ranks for testing multivariate elliptic white noise against ARMA dependence. *Bernoulli* **8** 787–815. MR1963662

Hallin, M. and Paindaveine, D. (2002c). Multivariate signed ranks: Randles' interdirections or Tyler's angles? In *Statistical Data Analysis Based on the $L_1$ Norm and Related Methods* (Y. Dodge, ed.) 271–282. Birkhäuser, Basel. MR2001322

Hallin, M. and Paindaveine, D. (2003). Affine invariant linear hypotheses for the multivariate general linear model with ARMA error terms. In *Mathematical Statistics and Applications: Festschrift for Constance van Eeden* (M. Moore, S. Froda and C. Léger, eds.) 417–434. IMS, Beachwood, OH.

Hallin, M. and Paindaveine, D. (2004a). Asymptotic linearity of serial and nonserial multivariate signed rank statistics. *J. Statist. Plann. Inference*. To appear.

Hallin, M. and Paindaveine, D. (2004b). Multivariate signed rank tests in vector autoregressive order identification. *Statist. Sci.* To appear. MR2001322

Hallin, M. and Paindaveine, D. (2005). Affine-invariant aligned rank tests for the multivariate general linear model with VARMA errors. *J. Multivariate Anal.* **93** 122–163.

Hallin, M. and Puri, M. L. (1988). Optimal rank-based procedures for time-series analysis: Testing an ARMA model against other ARMA models. *Ann. Statist.* **16** 402–432. MR924878

Hallin, M. and Puri, M. L. (1991). Time-series analysis via rank-order theory: Signed-rank tests for ARMA models. *J. Multivariate Anal.* **39** 1–29. MR1128669

Hallin, M. and Puri, M. L. (1994). Aligned rank tests for linear models with autocorrelated error terms. *J. Multivariate Anal.* **50** 175–237. MR1293044

Hallin, M. and Puri, M. L. (1995). A multivariate Wald–Wolfowitz rank test against serial dependence. *Canad. J. Statist.* **23** 55–65. MR1340961

Hallin, M. and Tribel, O. (2000). The efficiency of some nonparametric rank-based competitors to correlogram methods. In *Game Theory, Optimal Stopping, Probability and Statistics. Papers in Honor of Thomas S. Ferguson on the Occasion of His 70th Birthday* (F. T. Bruss and L. Le Cam, eds.) 249–262. IMS, Beachwood, OH. MR1833863

Hallin, M. and Werker, B. J. M. (1999). Optimal testing for semi-parametric AR models: From Gaussian Lagrange multipliers to autoregression rank scores and adaptive tests. In *Asymptotics, Nonparametrics and Time Series* (S. Ghosh, ed.) 295–350. Dekker, New York. MR1724702

Hallin, M. and Werker, B. J. M. (2003). Semiparametric efficiency, distribution-freeness and invariance. *Bernoulli* **9** 167–182. MR1963675

Hannan, E. J. (1970). *Multiple Time Series*. Wiley, New York. MR279952

Hettmansperger, T. P., Nyblom, J. and Oja, H. (1994). Affine invariant multivariate one-sample sign tests. *J. Roy. Statist. Soc. Ser. B* **56** 221–234. MR1257809

Hettmansperger, T. P., Möttönen, J. and Oja, H. (1997). Affine invariant multivariate one-sample signed-rank tests. *J. Amer. Statist. Assoc.* **92** 1591–1600. MR1615268

Hodges, J. L., Jr. and Lehmann, E. L. (1956). The efficiency of some nonparametric competitors of the $t$-test. *Ann. Math. Statist.* **27** 324–335. MR79383

Jan, S.-L. and Randles, R. H. (1994). A multivariate signed sum test for the one-sample location problem. *J. Nonparametr. Statist.* **4** 49–63. MR1366363

Le Cam, L. (1986). *Asymptotic Methods in Statistical Decision Theory*. Springer, New York. MR856411

Liebscher, E. (2005). A semiparametric density estimator based on elliptical distributions. *J. Multivariate Anal.* **92** 205–225. MR2102252

Möttönen, J. and Oja, H. (1995). Multivariate spatial sign and rank methods. *J. Nonparametr. Statist.* **5** 201–213. MR1346895

Möttönen, J., Oja, H. and Tienari, J. (1997). On the efficiency of multivariate spatial sign and rank tests. *Ann. Statist.* **25** 542–552. MR1439313

Möttönen, J., Hettmansperger, T. P., Oja, H. and Tienari, J. (1998). On the efficiency of affine invariant multivariate rank tests. *J. Multivariate Anal.* **66** 118–132. MR1648529

Oja, H. (1999). Affine invariant multivariate sign and rank tests and corresponding estimates: A review. *Scand. J. Statist.* **26** 319–343. MR1712063

Peters, D. and Randles, R. H. (1990). A multivariate signed-rank test for the one-sample location problem. *J. Amer. Statist. Assoc.* **85** 552–557. MR1141757

Puri, M. L. and Sen, P. K. (1971). *Nonparametric Methods in Multivariate Analysis*. Wiley, New York. MR298844

Randles, R. H. (1989). A distribution-free multivariate sign test based on interdirections. *J. Amer. Statist. Assoc.* **84** 1045–1050. MR1134492

Randles, R. H. (2000). A simpler, affine-invariant, multivariate, distribution-free sign test. *J. Amer. Statist. Assoc.* **95** 1263–1268. MR1792189

Rao, C. R. and Mitra, S. K. (1971). *Generalized Inverses of Matrices and Its Applications*. Wiley, New York. MR338013

Tyler, D. E. (1987). A distribution-free $M$-estimator of multivariate scatter. *Ann. Statist.* **15** 234–251. MR885734

Um, Y. and Randles, R. H. (1998). Nonparametric tests for the multivariate multi-sample location problem. *Statist. Sinica* **8** 801–812. MR1651509

Institut de Statistique et
de Recherche Opérationnelle
Université Libre de Bruxelles
Campus de la Plaine CP 210
B-1050 Bruxelles
Belgium
e-mail: mhallin@ulb.ac.be
e-mail: dpaindav@ulb.ac.be